\documentclass[11pt,leqno]{article}

\usepackage{ShVa787}

\begin{document}


%
\title{\TitleText}

\author{%
Saharon Shelah %
\thanks{\ShThanks}\\%
\and
Pauli V\"{a}is\"{a}nen %
\thanks{\PVThanks}\\%
}

\date{\today}
\maketitle

\begin{abstract}%
%
%
We strengthen non-structure theorems for almost free Abelian groups by
studying long Ehrenfeucht-Fra\"\i ss\'e games between a fixed group of
cardinality $\lambda$ and a free Abelian group.
A group is called $\epsilon$-game-free if the isomorphism player has a
winning strategy in the game (of the described form) of length
$\epsilon \in \lambda$.
We prove for a large set of successor cardinals $\lambda = \mu^+$ the
existence of nonfree $(\mu \cdot \omega_1)$-game-free groups of
cardinality \lambda. We concentrate on successors of singular
cardinals.
%
%
%
 \footnote{%
2000 Mathematics Subject Classification: %
\SubjectClassification
 Key words: %
\KeyWords
}%
\end{abstract}

%



\begin{SECTION}{Introduction}{Introduction}

%
The problem of possible cardinals carrying a nonfree almost free
Abelian group has already a long history, see e.g. \cite
{EklofMekler}.
Recall that if \Ga is a group of cardinality \lambda, \Ga is \InfLan
\lambda-equivalent to a free Abelian group if and only if \Ga is
strongly \lambda-free \cite {EklofInfLan}.
Recall also that the \InfLan \lambda-equivalence is characterized by
an \EF game of length \omega \Note {$< \lambda$ elements at the time}.
For an ordinal \epsilon, a group \Ga is called \epsilon-game-free, if
the ``isomorphism player'' has a winning strategy in the \EF game of
length \epsilon between \Ga and a free Abelian group \Note {countable
many elements at the time}.
Let \mu be a cardinal.  We study existence of nonfree groups of
cardinality \SuccCard \mu which are \epsilon-game-free with $\mu \leq
\epsilon < \SuccCard \mu$. We concentrate our attention on the case
that \mu is singular. This work continues \cite {PV1}, where the case
that \lambda is a successors of a regular cardinal is studied.
The following result is the main theorem of the paper \Note {presented
  in a more general form in \Section {Conclusion}}.

\begin{THEOREM}{SimplifiedForm}
Let \MgShSet denote the smallest set of cardinals such that \Alephf 0
is in \MgShSet and \MgShSet is closed under the operations $\lambda
\mapsto \SuccCard \lambda$ and $\Pair \lambda \kappa \mapsto \SuccCard
[\kappa +1] \lambda$.
For every cardinal \lambda of the form \SuccCard \mu in \MgShSet with
$\mu \geq \Alephf 2$, there exists a nonfree $(\mu \Times \Omegaf
1)$-game-free group of cardinality \lambda.
\end{THEOREM}

For an ordinal \epsilon and \lambda of the form \SuccCard \mu, $(\mu
\Times \epsilon)$-game-freeness \Note {the ordinal multiplication} of
a group \Ga implies that \Ga is equivalent to a free Abelian group
with respect to a ``deep'' infinitary language \InfLan [\theta]
\lambda, even a stronger language than \InfLan \lambda \cite
{Hytt90,Kartt,Oikk97}.
So our results can be interpreted as strengthened non-structure
theorems for almost free Abelian groups.

Shelah proved in \cite {Sh161} that the question of existence of
nonfree almost free Abelian groups is equivalent to a purely set
theoretical question concerning existence of nonfree almost free
families of countable sets. A family is almost free if all the
subfamilies of cardinality strictly less than the whole family has a
transversal. A transversal for a subfamily is an injective choice
function whose domain is the subfamily.
It is nowadays a standard custom to write that \NPT \lambda \kappa
holds
if there exists a family \SFamily of sets such that the elements of
\SFamily have cardinality $\leq \kappa$, \SFamily is almost free,
\SFamily is not free, and \SFamily has cardinality \lambda
\Note {in some papers \NPT \lambda {\SuccCard \kappa} is used instead
of \NPT \lambda \kappa}.
For some history of \NPT \lambda \kappa see, e.g., \cite [II.\SS0]
{CardArith} and \cite [\SS0] {MgSh204}.
Recall that \NPT \lambda \kappa fail for all $\lambda > \kappa$ if
\lambda is a singular cardinal \cite {Sh52}.
Recall also that by \cite [\SS1] {Sh521}, for every cardinal \lambda,
\NPT \lambda {-} holds iff there exists a nonfree, \lambda-separable
group of cardinality \lambda.

A transversal game on a family \SFamily of sets is a two players game,
where on each round $i$ the first player chooses a set $s_i$ from
\SFamily and the second player, also called transversal player, must
answer with an element $x_i$ from $s_i$ so that $x_i$ is distinct from
the elements $x_j$, $j < i$, the second player has chosen on the
earlier rounds \Note {\Definition {game-free_SFamily}}.  That means
that the transversal player must be able to choose an ``extendable''
transversal whose domain contains at least the sets chosen by the
first player.
A family \SFamily is called \epsilon-game-free if the transversal
player has a winning strategy in the transversal game of length
\epsilon on \SFamily.
As can be expected, the question
``does there exist a nonfree \epsilon-game-free group of cardinality
\lambda'',
for various \epsilon,
is very closely related to the question ``does there exist a nonfree
\epsilon-game-free family of countable sets having cardinality
\lambda''.
For a detailed exposition of a transformation of an \epsilon-game-free
family into an \epsilon-game-free group see \cite [\AuxRef {PV1}
{Subsection} {StronglyFree}] {PV1}.

Let \chi denote the first cardinal fixed point \Note {the first
cardinal \kappa with $\Alephf \kappa = \kappa$}.
Our main target is to prove, without any assumptions beyond \ZFC, that
for a given $\theta < \chi$, there are nonfree $\mu \Times
\theta$-game-free groups of cardinality \SuccCard \mu for unbounded
many singular cardinals \mu below \chi.
In \cite [\SS1] {MgSh204} it is proved that \NPT {\SuccCard \mu} {-}
holds for unbounded many singular cardinals \mu below \chi.
The main theorem follows from this result together with a proposition
on ``canonical'' families of countable sets \Note {\Section
{TransGame}}.

We present definitions of \lambda-sets and \lambda-systems in
\Section {Preliminaries} even though those could be found from \cite
[\SS3] {Sh161}, \cite [Appendix] {Sh521}, or \cite [VII.3A]
{EklofMekler}. The reason for this is that we want the new terminology
``canonical form of a \lambda-set'' and ``\NPT \lambda {-}-skeletons
of type \seq \lambda'' \Note {or briefly ``incompactness skeletons''}
to be very clear for the reader. Moreover, we introduce a variant of
these, called ``\NRT \lambda {-}-skeletons of type \seq \lambda'', in
\Section {NRT}.

In \Section {TransGame} we prove the main proposition. Namely we show
that any family \SFamily of countable sets of cardinality \SuccCard
\mu with \mu a singular cardinal, whose canonical form fulfills a
certain ``cofinality'' condition, can be transformed into a family
\SFamily['] of countable sets such that \SFamily['] is a nonfree
\mu-game-free family having cardinality \SuccCard \mu.
The reader may wonder why the conclusion of \Proposition {Singular} is
\epsilon-game-free for every $\epsilon < \mu$ instead of
\mu-game-free. The answer is that by \cite [\AuxRef {PV1} {Lemma}
{SingularCompactness}] {PV1}, \epsilon-game-freeness for every
$\epsilon < \mu$ implies \mu-game-freeness, when \mu is a singular
cardinal.

By \cite [\SS3] {Sh521} it is consistent, relative to existence of two
weakly compact cardinals, that for some uncountable regular cardinals
$\kappa < \lambda$, both \NPT \lambda \kappa and \NPT \kappa {-} hold,
even though, \NPT \lambda {-} does not hold. Hence the ordinary notion
of ``almost free'' is not strong enough for the transitivity
conclusion: for all $\lambda_1 > \lambda_2 > \lambda_3$, if \NPT
{\lambda_1} {\lambda_2} and \NPT {\lambda_2} {\lambda_3} hold, then
\NPT {\lambda_1} {\lambda_3} hold.
In \Section {NRT} we present special kind of families, called \NRT -
{\kappa}-families, for which an analogical transitivity property hold
\Note {$R$ refers to the notion \NRT \lambda \kappa, which will be an
analog of \NPT \lambda \kappa}: for all $\lambda_1 > \lambda_2 >
\lambda_3$, \NRT {\lambda_1} {\lambda_2} and \NRT {\lambda_2}
{\lambda_3} implies \NRT {\lambda_1} {\lambda_3}.
Note that if \SFamily exemplifies \NPT \lambda \kappa, and \SFamily
satisfies a stronger freeness notion, say the following property:
\begin{property}

whenever \R is a subset of \SFamily of cardinality $< \lambda$, there
are pairwise disjoint sets \Seq {r_s} {s \in \R} such that $s \Minus
r_s$ has cardinality $< \kappa$ for every $s \in \R$;

\end{property}
then \SFamily and any family exemplifying \NPT \kappa {-} can be
amalgamated to get a family exemplifying \NPT \lambda {-}.
The point of introducing \NRT - \kappa-families and ``\NRT [\Ideal] -
\kappa-freeness'' in \Section {NRT} is that the existence of such a
family is also necessary to build a canonical example of a family of
countable sets.
In other words, the largest nicely incompact subset of \chi \Note
{which denotes the first cardinal fixed point}, defined in \cite
[\SS2] {Sh521}, coincides with the smallest set of cardinals below
\chi which contains \Alephf 0 and is closed under ``amalgamation of
\NRT - \kappa-families'' \Note {\Definition {NRTSet}}.

The main theorem of the paper is presented in \Section {Conclusion}.
Namely we link together the pieces proved in \cite [\AuxRef {PV1}
{Section} {Transversals}] {PV1}, \Section {TransGame}, and
\Section {NRT}.
For all nonzero $n < \omega$, it is possible to build nonfree
$(\omega_n \Times \omega_{n-1})$-game-free families of cardinality
\Alephf {n+1} level by level using the ``old methods'', see \cite
[\AuxRef {PV1} {Section} {Basic}] {PV1}.
To get an example of nonfree \mu-game-free group of cardinality
\SuccCard \mu with $\mu = \Alephf \omega$, one needs \Proposition
{Singular} and \cite [Theorem 4 of \SS1] {MgSh204}.
To get examples for $\Alephf \omega < \mu < \Alephf {\omega +
  \omega}$, one may use an example for \Alephf \omega and apply \cite
[\AuxRef {PV1} {Lemma} {RegForm}] {PV1}, which just says that a
certain type of an \NPT {\SuccCard {\mu_2}} {-}-skeleton is $(\mu_2
\Times \mu_1)$-game-free, if $\mu_2 = \SuccCard {\mu_1}$ and the
``previous level'' is $\mu_1$-game-free
\Note {a group whose \Gamma-invariant consists of ordinals of
  cofinality $\geq \mu_2$ is always $\mu_2$-game-free, but not
  necessarily $\mu_2 + \mu_2$-game-free \cite [\AuxRef {PV1} {Section}
  {Basic}] {PV1}}.
By \cite [Theorem 4 of \SS1] {MgSh204} and \Proposition {Singular},
one get examples for even larger \mu's, e.g., $\mu = \Alephf
{\omega_n}$ for any $n < \omega$. In fact it follows that such
examples exist for unbounded many \mu's below the first cardinal fixed
point.

Suppose that \Ga is an Abelian group of cardinality \SuccCard \mu
whose \Gamma-invariant \GammaInv \Ga contains stationary many \theta
cofinal ordinals, i.e., for some filtration \Seq {\Gaf \alpha} {\alpha
  < \SuccCard \mu} of \Ga, the set $\Cof \theta \Inter \Set {\alpha <
  \SuccCard \mu} {\Quotient {\Gaf {\alpha+1}} {\Gaf \alpha} \Text {is
    not free}}$ is stationary in \SuccCard \mu.
Assume, for a while, that $\GammaInv \Ga \Inter \Cof \theta$ is in the
ideal \IGoodIdeal {\SuccCard \mu} of ``all good subsets of \SuccCard
\mu''
\Note {see, e.g., \cite [Analytical Guide \SS0] {CardArith}}.
Then by \cite [\AuxRef {PV1} {Lemma} {GameFreeImpliesFree}] {PV1},
$(\mu \Times (\theta + 1))$-game-freeness implies freeness for groups
\Ga of cardinality \SuccCard \mu.
Hence \Theorem {game-free} is very close to the optimal result
provable from \ZFC alone.

What happens if \mu is singular and \GammaInv \Ga is not in
\IGoodIdeal {\SuccCard \mu}? \Note {For the regular case see \cite
  [\AuxRef {PV1} {Proposition} {CubGame}] {PV1}.} In \Section
{CubGame} we collect together some facts about ``good points with
respect to a scale for \mu'' and we prove that, relative to the
existence of a supercompact cardinal,
it is possible to have a nonfree group \Ga of cardinality \SuccCard
\mu \Note {e.g. $\mu = \Alephf {\omega_n}$ for any $n < \omega$}, so
that \Ga is \epsilon-game-free for every $\epsilon < \SuccCard
\mu$.
It is also possible to obtain the ``maximal game-freeness'' without
\Ga being free by measuring the length of the game by trees.
%

\end{SECTION}


\begin{SECTION}{Preliminaries}{Preliminaries}

%
We denote the cofinality of an ordinal \alpha by \Cf \alpha and the
cardinality of a set $X$ by \Card X. The class of all ordinals of
cofinality \theta is denoted by \Cof \theta. The set of all subsets of
$X$ of cardinality $< \lambda$ is denoted by \SubsetsOfCard X {<
\lambda}.
%

The definitions concerning \lambda-sets and \lambda-systems are from
\cite [\SS3] {Sh161}. There are slightly revised versions in \cite
[Appendix] {Sh521}. For the reader's convenience, we have chosen the
notation so that it is compatible to \cite [\SS3 of Chapter VII]
{EklofMekler}, and of course, compatible to \cite {PV1}.
The reason to represent the definitions is that in the later sections
the concepts of ``type of a \lambda-set'' and ``incompactness
skeletons'' have a central role.

\begin{DEFINITION}{Almost_free_family}%
%
%
Suppose \SFamily is a family of countable sets. A function $T$ from
\SFamily into \BigUnion \SFamily is called a transversal for \SFamily
when $T$ is injective and for every $s \in \SFamily$, $T(s) \in s$.
If \SFamily is enumerated by \Set {\s i} {i \in I} without repetition,
then a transversal for $J \Subset I$ means a transversal $T$ for \Set
{\s i} {i \in J}. When more convenient, we abbreviate $T(\s i)$ by
$T(i)$.

A family \SFamily is called free if there exists a transversal for
\SFamily. For a cardinal \kappa, \SFamily is called \kappa-free, when
every subfamily of \SFamily of cardinality $< \kappa$ is
free. \SFamily is almost free when it is \lambda-free with $\lambda =
\Card \SFamily$.

For a subfamily \R of \SFamily, \Quotient \SFamily \R denotes the
family \Set {s \Minus \Union \R} {s \in \SFamily \Minus \R}.
\end{DEFINITION}

\begin{DEFINITION}{game-free_SFamily}%
Suppose \SFamily is a family of countable sets and \epsilon is an
ordinal. We call \SFamily \epsilon-game-free if the second player,
called \PlayerII, has winning strategy in the following two players
``transversal game'' denoted by \TransGame \epsilon \SFamily.
A play of the game last at most \epsilon rounds. On each round $i$ the
first player, called \PlayerI, chooses a countable subfamily $\R_i$ of
\SFamily.
The second player, \PlayerII, must answer with a transversal $T_i$
whose domain contains all the elements of $\R_i$ and which extends the
transversals $T_j$, $j < i$, he has chosen on the earlier
rounds. \PlayerII [big] wins a play if he succeeds to follow the given
rules \epsilon rounds.
\end{DEFINITION}

\begin{DEFINITION}{LambdaSet}%
\AddInTheoremCite[Definition 3.1]{Sh161}%
Suppose \lambda is an uncountable regular cardinal, and \LS is a
nonempty set of finite sequences of ordinals which is closed under
initial segments.
Let \Sf denote the ``{\,}final nodes'' of \LS, i.e., \Sf is the
smallest subset of \LS with $\LS = \Set [\big] {\Res \eta m} {\eta \in
\Sf \And m \leq \Lh \eta}$.
A set \LS as above is called a \lambda-set if there exist uncountable
regular cardinals \Seq {\lambda_\rho} {\rho \in \LS \Minus \Sf} such
that
\begin{properties}

\item $\lambda_\emptyset = \lambda$;

\item for every $\rho \in \LS \Minus \Sf$,
$
        \LSE \rho
        =
        \Set {\alpha}
        {\ConcatSimple \rho \alpha \in \LS}
$
is a stationary subset of $\lambda_\rho$ \Note {so $\ConcatSimple \rho
\alpha \in \LS$ implies that $\alpha < \lambda_\rho$};

\item for every $\rho \in \LS \Minus \Sf$ and $\alpha \in \LSE \rho$,
$\lambda_{\ConcatSimple \rho \alpha} \leq \Card \alpha$ \Note {so
$\lambda_\rho > \lambda_{\ConcatSimple \rho \alpha}$}.

\end{properties}
The sequence \Seq {\lambda_\rho} {\rho \in \LS \Minus \Sf} is called
the type of \LS.
A \lambda-set \LS is said to have height $n$ if $n$ is a finite
ordinal such that $\Lh \eta = n$ for all $\eta \in \Sf$.
\end{DEFINITION}

\begin{DEFINITION}{Small}%
Suppose \LS is a \lambda-set of type \Seq {\lambda_\rho} {\rho \in \LS
\Minus \Sf}. 
A subset \LS['] of \LS is called a sub-\lambda-set of \LS if \LS['] is
a \lambda-set such that $\Sf['] \Subset \Sf$ and the type of \LS['] is
the restriction \Seq {\lambda_\rho} {\rho \in \LS['] \Minus \Sf[']} of
the type of \LS.
A subset $I$ of \Sf is small in \Sf if the set \Set {\Res \eta m}
{\eta \in I \And m \leq \Lh \eta} is not a sub-\lambda-set of \LS.
\end{DEFINITION}

\begin{DEFINITION}{CanonicalLambdaSet}%
\AddInTheoremCite[Claim 3.2]{Sh161}%
We say that a \lambda-set \LS of type \seq \lambda has a canonical
form when the following demands are fulfilled.
\LS has height \Star n and
there exist sequences of cardinals \Seq {\lambda_n'} {n < \Star n} and
\Seq {\theta_n} {n < \Star n} such that
$\lambda_0' = \lambda_\emptyset$
and for every $\rho \in \LS \Minus \Sf$ of length $n$,
\begin{enumprops}

\ITEM{RegLimit}%
either both of the following two properties are satisfied:
\begin{properties}

\item the set \LSE \rho consist of regular limit cardinals and
$\theta_n = 0$ \Note {in this case $\theta_n$ is called undefined};

\item if $n+1 < \Star n$ then for every $\alpha \in \LSE \rho$,
$\lambda_{\ConcatSimple \rho \alpha} = \alpha$ and $\lambda_{n+1} = 0$
\Note {in this case $\lambda_{n+1}'$ is called undefined};

\end{properties}

\ITEM{Otherwise}%
or otherwise both of the following two properties are satisfied:
\begin{properties}

\item \LSE \rho is a subset of \Cof {\theta_n};

\item if $n+1 < \Star n$ then for every $\alpha \in \LSE \rho$,
$\lambda_{\ConcatSimple \rho \alpha} = \lambda_{n+1}'$.

\end{properties}

\end{enumprops}
\end{DEFINITION}

\begin{FACT}{CanonicalLambdaSet}%
\AddInTheoremCite[Claim 3.2]{Sh161}%
Every \lambda-set contains a sub-\lambda-set which is in a canonical
form.
\end{FACT}

\begin{DEFINITION}{LambdaSystem}%
\AddInTheoremCite[Definition 3.4]{Sh161}%
Suppose \lambda is an uncountable regular cardinal, \LS is a
\lambda-set of type $\seq \lambda = \Seq {\lambda_\rho} {\rho \in \LS
\Minus \Sf}$.
An indexed family
\[
        \LSystem =
        \Seq [\big] {\LSet {\Concat \rho {\SimpleSeq \alpha}}}
        {\rho \in \LS \Minus \Sf \And \alpha < \lambda_\rho}
\]
is called a \lambda-system of type \seq \lambda, if \LSet \emptyset is
the empty set
\Note {only for technical reasons}
and for every $\rho \in \LS \Minus \Sf$ the sequence
$
        \seq {\LSet \rho}
        =
        \Seq [\big]
        {\LSet {\ConcatSimple \rho \alpha}}
        {\alpha < \lambda_\rho}
$
satisfies that
\begin{properties}

\item \seq {\LSet \rho} is a strictly increasing continuous chain of
sets;

\item the union of \seq {\LSet \rho} has cardinality $\lambda_\rho$;

\item each \LSet {\ConcatSimple \rho \alpha} has cardinality $<
\lambda_\rho$.


\end{properties}
A \lambda-system is called disjoint when the sets \BigUnion [\alpha <
\lambda_\rho] {\LSet {\ConcatSimple \rho \alpha}}, for all $\rho \in
\LS \Minus \Sf$, are disjoint.
\end{DEFINITION}

\begin{DEFINITION}{NPTSkeleton}%
Suppose $\lambda > \kappa$ are infinite regular cardinals.
A tuple \NPTSkeleton is called an \NPT \lambda \kappa-skeleton \Note
{of type \seq \lambda and of height \Star n} when
\begin{properties}

\item \LS is a \lambda-set of type \seq \lambda;

\item \LS is in a canonical form and its height is \Star n;

\item the cardinals in \seq \lambda are greater than \kappa;

\item $\LSystem = \Seq [\big] {\LSet {\ConcatSimple \rho \alpha}}
{\rho \in \LS \Minus \Sf \And \alpha < \lambda_\rho}$ is a disjoint
\lambda-system of type \seq \lambda;

\item \SFamily is a family of sets of cardinality $\leq \kappa$
enumerated by \Set {\s \eta} {\eta \in \Sf};

\item \SFamily is based on \LSystem, which means that for every $\eta
\in \Sf$,
\[
        \s \eta \Subset
        \BigUnion
        [m \leq \Star n]
        {\LSet {\Res \eta m}}.
\]

\end{properties}
Additionally, the demands in \cite [Claim 3.6 and 3.7] {Sh161} are
fulfilled
\Note {in case $\kappa = \Alephf 0$ they are essentially equivalent to
the ``beautifulness properties'' presented in \cite [Definition
VII.3A.2(1--6)] {EklofMekler}}.
Define for every $\eta \in \Sf$ and $m < \Star n$, that
\[
        \slevel \eta {m+1} = \s \eta \Inter \LSet {\Res \eta {m+1}}.
\]
\begin{REMARK}
\slevel \eta {m+1} is denoted by \slevel \eta m in \cite {Sh161}.
\end{REMARK}

For the reader's convenience we repeat three of the additional demands
\Note {used in \Section {TransGame}\,}:
\begin{textprops}

\TextItem{det_by_index}
For all $\eta \not= \nu \in \Sf$, if $\s \eta \Inter \s \nu \not=
\emptyset$ then there is an unique $m < \Star n$ such that $\s \eta
\Inter \s \nu = \slevel \eta {m+1} \Inter \slevel \nu {m+1}$ and for
every $l < \Star n$ with $l \not= m$, $\eta(l) = \nu(l)$.

\TextItem{UncountableCof}
For every $\eta \in \Sf$ and $m < \Star n$, if $\eta(m)$ has
cofinality $\theta > \kappa$, then there is $n < \Star n$ with
$\lambda_{\Res \eta n} = \theta$
\Note {$\lambda_{\Res \eta n}$ is from the type \seq \lambda of \LS}.

\TextItem{tree-order}
In case $\kappa = \Alephf 0$:
For all $\eta \in \Sf$ and $m < \Star n$, the sets \slevel \eta {m+1}
have enumerations \Seq {x^{\eta,m}_l} {l < \omega} with the following
property:
for every $\nu \in \Sf$ and $n < \omega$,
if $x^{\nu,m}_n \in \slevel \eta {m+1}$ then
the initial segments \Seq {x^{\nu,m}_l} {l \leq n} and \Seq
{x^{\eta,m}_l} {l \leq n} are equal.

\end{textprops}
\end{DEFINITION}

\begin{DEFINITION}{scale}%
Suppose $\seq \mu = \Seq {\mu_\xi} {\xi < \kappa}$ is an increasing
sequence of regular cardinals approaching \mu such that $\mu_0 >
\SuccCard \kappa$.
For functions $f, g \in \Prod [\xi < \kappa] {\mu_\xi}$, we write that
$f \EventLess g$ when \Def {$f$ is eventually strictly less than $g$},
i.e., if \Set [\big] {\xi < \kappa} {f(\xi) \geq f(\xi)} is bounded in
\kappa.

A pair \Pair {\seq \mu} {\scale f} is called a scale for \mu when \seq
\mu is as above and $\scale f = \Seq {f_\alpha} {\alpha < \lambda}$ is
a \EventLess-increasing cofinal sequence of functions in \Prod [\xi <
\kappa] {\mu_\xi}
\Note {cofinal means that for every $g \in \Prod [\xi < \kappa]
{\mu_i}$ there exists \alpha with $g \EventLess f_\alpha$}.
\end{DEFINITION}

\begin{FACT}{scale}%
\AddInTheoremCite [Theorem II.1.5] {CardArith}%
For any singular cardinal \mu, there exists a scale for \mu.
In other words, there is a sequence \Seq {\mu_\xi} {\xi < \kappa} of
regular cardinals approaching \mu so that \Prod [\xi < \kappa] {\mu_i}
has true cofinality \SuccCard \mu under the partial order \EventLess
\Note {denoted by $<_{J^{\Symb{bd}}_k}$ in the reference}.
\end{FACT}
%

\end{SECTION}


\begin{SECTION}{TransGame}{Transversals and game-free families}

%
The paper \cite {PV1} left open a question: if \mu is a singular
cardinal, does there exist a family \SFamily of countable sets of
cardinality \SuccCard \mu which is nonfree and \epsilon-game-free for
all $\epsilon < \mu$?  Shelah has a very satisfactory answer to this
question.

\begin{PROPOSITION}{Singular}%
%
%
Suppose \NPTSkeleton ['] is an \NPT \lambda {-}-skeleton of type $\seq
\lambda = \Seq {\lambda_\rho} {\rho \in \LS['] \Minus \Sf[']}$ such
that the following conditions hold:
\begin{properties}

\item%
\mu is a singular cardinal and \lambda equals \SuccCard \mu;

\item%
either $\Cf \mu = \Alephf 0$
or there is $m$ below the height of \LS such that $\lambda_\rho = \Cf
\mu$ for every $\rho \in \LS$ of length $m$;

\item%
the cardinal \theta for which $\LSE ['] \emptyset \Subset \Cof \theta$
holds \Note {and which exists by \Definition {NPTSkeleton}} is such
that $\Cf \mu \not= \theta$.

\end{properties}
Then there exists \NPT \lambda {-}-skeleton \NPTSkeleton such that
\SFamily is \epsilon-game-free for every $\epsilon < \mu$.
Moreover, \LS is a sub-\lambda-set of \LS['] and hence $\LSE \emptyset
\Subset \Cof \theta$ holds \Note {$\LS = \LS[']$ if \seq \lambda does
not contain limit cardinals}.
\end{PROPOSITION}

The rest of the section is devoted to the proof of this proposition.

Suppose that \LS['] has height \Star n \Note {by the canonical form
  height is well-defined}. Now note that \Star n must be at least $2$:
If \mu has uncountable cofinality, this follows directly from our
second demand.
On the other hand, if \mu has countable cofinality, the claim follows
from the demand that $\LSE['] \emptyset \not\Subset \Cof {\Alephf 0}$
together with \RefOfDefinition {NPTSkeleton} {UncountableCof}.

\begin{REMARK}
If $\Star n = 1$ and $\LSE \emptyset \Subset \Cof {\Alephf 0}$, then
\SFamily would not be \Alephf 1-game-free as explained in \cite
[\AuxRef {PV1} {Example} {not_strongly_game-free}] {PV1}.
\end{REMARK}

We let $\lambda_1$ denote the regular cardinal given by \Fact
{CanonicalLambdaSet}, i.e., for every $\alpha \in \LSE \emptyset$,
$\lambda_{\SimpleSeq \alpha} = \lambda_1$. By taking a suitable
sub-\lambda-set \LS of \LS['] if necessary, we may assume that
\begin{property}

if \Cf \mu is uncountable, then there is a fixed $m < \Star n$ such
that for every $\eta \in \Sf$, $\lambda_{\Res \eta m} = \Cf \mu$.

\end{property}
If there is no limit cardinals in \seq \lambda, then \Fact
{CanonicalLambdaSet} guarantees that we can choose $\LS = \LS[']$.

By \Fact {scale} choose a scale \Pair {\seq \mu} {\scale f} for \mu
\Note {\Definition {scale}}, where $\seq \mu = \Seq [\big] {\mu_\xi}
{\xi < \Cf \mu}$ and $\scale f = \Seq {f_\alpha} {\alpha < \lambda}$.

The only modification of the family $\SFamily ['] = \Set {\s \eta'}
{\eta \in \Sf [']}$ needed is that the ``first coordinate'' $(\slevel
\eta 1)' = \s \eta' \Inter \LSet ['] {\Res \eta 1}$ of each $\s \eta'$
is slightly changed. The modification depends on the cofinality of
\mu as follows:
Fix \eta from \Sf. If \mu has countable cofinality, then define
\[
        \slevel \eta 1 = 
        \Set {\Res {F_\eta} l} {l < \omega},
\]
where $F_\eta$ is a fixed function from \omega onto $\Ran
{f_{\eta(0)}} \times (\slevel \eta 1)'$. This definition ensures that
the property \RefOfDefinition {NPTSkeleton} {tree-order} becomes
fulfilled.
In order to choose the corresponding new \lambda-system, let \LSet
{\SimpleSeq \alpha}, for every $\alpha < \lambda$, be $\SubsetsOfCard
[\big] {\mu \times \LSet ['] {\SimpleSeq \alpha}} {< \Alephf 0}$,
and additionally, for all nonempty $\rho \in \LS \Minus \Sf$ and
$\alpha < \lambda_\rho$, define \LSet {\ConcatSimple \rho \alpha} to
be the old \LSet ['] {\ConcatSimple \rho \alpha}.
Otherwise, for the fixed \eta, $\Cf \mu = \lambda_{\Res \eta m} >
\Alephf 0$ holds and we define
\[
        \slevel \eta 1 = 
        \Singleton [\Big] {f_{\eta(0)}\Par[\big]{\eta(m)}}
        \times (\slevel \eta 1)',
\]
and for every $\alpha < \lambda$, $\LSet {\SimpleSeq \alpha} = \mu
\times \LSet ['] {\SimpleSeq \alpha}$
\Note {\RefOfDefinition {NPTSkeleton} {tree-order} is satisfied since
$(\slevel \eta 1)'$ satisfies it}.
Note that now \NPTSkeleton is an \NPT \lambda {-}-skeleton, i.e., it
fulfills all the demands mentioned in \Definition {NPTSkeleton},
because \NPTSkeleton ['] is an \NPT \lambda {-}-skeleton.

In order to show that \PlayerII has a winning strategy in the game
\TransGame \epsilon \SFamily for every $\epsilon < \mu$, it suffices
to describe a winning strategy for \PlayerII in a modified game of
length \sigma for every regular cardinal \sigma with $\lambda > \sigma
> \lambda_1$,
where the rules of the modified game are exactly as in \TransGame
\sigma \SFamily, except that \PlayerII is demanded to choose a
transversal only if the index of the round is a limit ordinal.

On every round $i < \sigma$ \PlayerII chooses two elementary submodels
\M {i+1} and \N {i+1} of \HerStr \chi {\seq f, \seq \mu, \LS,
  \LSystem, \SFamily} such that
\begin{properties}

\item \chi is some large enough regular cardinal;

\item $\Card {\M {i+1}} = \mu$;

\item $\Card {\N {i+1}} = \sigma$;

\item $\N {i+1} \Subset \M {i+1}$;

\item \N {i+1} contains all the elements chosen by \PlayerI;

\item $\sigma +1 \Subset \N {i+1}$;

\item $\M {i+1} \Inter \lambda$ is an ordinal denoted by $\delta_{i+1}
\in \lambda$;

\item $\M i, \N i \in \N {i+1}$, where $\M 0 = \N 0 = \emptyset$ and
for limit $i$, $\M i = \BigUnion [j < i] {\M j}$ and $\N i = \BigUnion
[j < i] {\N j}$.

\end{properties}

For the rest of this proof assume $i < \sigma$ to be a limit ordinal
or zero, and suppose \M j, \N j for each $j \leq i$ are chosen.
Denote the set $\SFamily \Inter \N i$ by \R.
We show that
\begin{mathprop}{subclaim}%
\Quotient \SFamily \R \Text{is \lambda-free},
\end{mathprop}%
because then after $i + \omega < \sigma$ rounds, \PlayerII is able to
continue \Note {or start} with some transversal $T_{i + \omega}$ whose
domain consist of the elements in $\SFamily \Inter \N {i + \omega}$
\Note {contains the elements chosen by \PlayerI so far}
and satisfies that $T_i \Subset T_{i+\omega}$ \Note {where $T_0 =
\emptyset$ and for $i$ which is a limit of limit ordinals, $T_i =
\BigUnion [j < i] {T_{j + \omega}}$}.

To check the details, we have to define some auxiliary notations \Note
{familiar from \cite {PV1} or \cite [VII.3A] {EklofMekler}}:
for every $\alpha < \lambda$, 
\[\begin{array}{lcl}
\S \alpha &=&
\Set [\big] {\s \eta} {\eta \in \Sf \And \eta(0) < \alpha};\\
\SfProj {\SimpleSeq \alpha} &=& \Set [\big] {\eta \in \Sf} {\eta(0) = \alpha};\\
\SFamilyProj {\SimpleSeq \alpha} &=&
\Set [\big] {\BigUnion [0 < l < \Star n] {\slevel \eta {l+1}}}
{\eta \in \SfProj \alpha}.\\
\end{array}\]
Hence, e.g., $\SFamily \Inter \M i$ equals to \S {\delta_i} for every
$i < \sigma$. We also need the fact that when \NPTSkeleton is an \NPT
\lambda {-}-skeleton, for every \alpha in \LSE \emptyset and for all
$I \Subset \SfProj {\SimpleSeq \alpha}$
\Note {recall smallness from \Definition {Small}}:
\begin{textprop}{Trans}%
$I$ is small in \SfProj {\SimpleSeq \alpha} iff there is a transversal
$T$ such that its domain is \Set {\s \eta} {\eta \in I} and for every
$\eta \in I$, $T(\s \eta) \not\in \slevel \eta 1$.
\end{textprop}%
Much more is proved in \cite [Claim 3.8] {Sh161}. This simple fact is
explained, e.g., in \cite [\AuxRef {PV1} {Fact} {TransProp}] {PV1}.

The proof of \Property {subclaim} is divided into several parts.
Because we may assume that \LSE \emptyset contains only limit
ordinals, \Quotient {\SFamily} {\S {(\delta_i) +1}} is
\lambda-free. Since $\R \Subset \S {\delta_i}$ it suffices to show
that
\Quotient {\S {(\delta_i) +1}} \R is free.
First we show that 
\begin{textprop}{above_delta_i}%
there is a transversal $T$ for \SFamilyProj {\SimpleSeq {\delta_i}} so
that $\Ran T \Inter \BigUnion \R$ is empty.
\end{textprop}%
\begin{REMARK}
We do not claim that $T$ witnesses \Quotient {\S {(\delta_i) +1}} {\S
  {\delta_i}} to be free. It can happen that \Ran T contains elements
from \BigUnion {\S {\delta_i}}.
\end{REMARK}

Secondly we prove that 
\begin{textprop}{below_delta_i}%
$T$ can be extended to a transversal $T'$ for \Quotient {\S
{(\delta_i) +1}} \R.
\end{textprop}

Note that if $\delta_i \not\in \LSE \emptyset$, then \Quotient {\S
{(\delta_i)+1}} {\S {\delta_i}} is empty, and \Property
{above_delta_i} holds trivially.
So suppose that $\delta_i \in \LSE \emptyset$. We want that
\begin{textprop}{f_delta_i}%
        there exists \Star \xi such that
        $f_{\delta_i}(\xi) \not\in \N i$
        for any $\xi \geq \Star \xi$.
\end{textprop}%
\relax From this the claim \Property {above_delta_i} follows in the following
way: By the definition of the ``new first coordinate'' there are two
different cases according to the cofinality of \mu.
If \mu has countable cofinality then choose a transversal $T'$ for
\Set [\big] {\s \eta'} {\eta \in \SfProj {\SimpleSeq {\delta_i}}}
\Note {$\s \eta'$ is the ``old first coordinate''}, and define the
desired transversal $T$ for \Set [\big] {\s \eta} {\eta \in \SfProj
{\SimpleSeq {\delta_i}}} by setting for every $\eta \in \SfProj
{\SimpleSeq {\delta_i}}$, that $T(\s \eta)$ is the finite restriction
of $F_\eta$ whose greatest element is the pair \Pair [\big]
{f_{\eta(0)}(\Star \xi)} {T'(\s \eta')}.
Suppose then that \Cf \mu is uncountable. By the definition of \slevel
\eta 1, when $\eta \in \SfProj {\SimpleSeq {\delta_i}}$, and the
assumption \Property {f_delta_i},
the following set is small in \SfProj {\SimpleSeq {\delta_i}}:
\[
        I_0 = \Set [\big] {\eta \in \SfProj {\SimpleSeq {\delta_i}}}
        {\slevel \eta 1 \Inter (\Union \R) \not= \emptyset}.
\]
Choose any transversal $U_1$ for $I_1 = \SfProj {\SimpleSeq
{\delta_i}} \Minus I_0$ \Note {$I_1 \Subset \Sf$ has cardinality $<
\lambda$}. By \Property {Trans}, $I_1' = \Set {\eta \in I_1}
{U_1(\eta) \not\in \slevel \eta 1}$ must be small in \SfProj
{\SimpleSeq {\delta_i}}.
Now the set $I_0 \Union I_1'$ is small in \SfProj {\SimpleSeq
{\delta_i}}, and by \Property {Trans} again, there is a transversal
$U_0$ for $I_0 \Union I_1'$ with the property that $U_0(\eta) \not\in
\slevel \eta 1$ for all $\eta \in I_0 \Union I_1'$.
Because of \RefOfDefinition {NPTSkeleton} {det_by_index}, the union $T
= U_0 \Union \Res {U_1} {I_1 \Minus I_1'}$ forms the desired
transversal witnessing that \Property {above_delta_i} holds. So
\Property {f_delta_i} implies \Property {above_delta_i}.


Note that \Set {\delta_j} {j < i} is a subset of \N i.
Define \Star \zeta to be the first index in \Cf \mu with $\mu_{\Star
  \zeta} > \sigma$.
For every $j \leq i$ define a function $h_j \in \Prod [\xi < \Cf \mu]
{\mu_\xi}$ by setting $h_j(\xi) = 0$ if $\xi < \Star \zeta$, and
$h_j(\xi) = \sup (\mu_\xi \Inter \N j)$ otherwise.
Clearly $h_k(\xi) < h_j(\xi)$ when $k < j \leq i$ and $\xi \geq \Star
\zeta$.
For every $j < i$ there is $\beta < \delta_i$ with $h_j \EventLess
f_\beta$, since $h_j \in \N i$ and $\N i \Inter \lambda \Subset
\delta_i$. Thus $h_j \EventLess f_{\delta_i}$ for every $j < i$.
To prove \Property {f_delta_i} it suffices to show that $h_i \EventLeq
f_{\delta_i}$.
For every $j < i$, there is the smallest $\xi_j < \Cf \mu \Minus \Star
\zeta$ satisfying that $\Set [\big] {\xi < \Cf \mu} {h_j(\xi) \geq
f_{\delta_i}(\xi)} \Subset \xi_j$.
By the definition of $h_j$'s, $\xi_k < \xi_j$ for every $k < j < i$.
By the assumption $\LSE ['] \emptyset \Inter \Cof {\Cf \mu} =
\emptyset$, $\Cf i = \Cf {\delta_i} \not= \Cf \mu$.
Hence there exists $\Star \xi \geq \Star \zeta$ such that $\xi_j \leq
\Star \xi$ for every $j < i$.
By the definition of $\xi_j$'s, for every $\xi \geq \Star \xi$,
$f_{\delta_i}(\xi) \geq \sup \Set {h_j(\xi)} {j < i} = h_i(\xi)$.
Therefore \Property {f_delta_i} holds, and we have proved \Property
{above_delta_i}.


Next we prove that there is some transversal for \Quotient {\S
{\delta_i}} \R, and finally we explain how to find a transversal
witnessing \Property {below_delta_i}.
In order to show that \Quotient {\S {\delta_i}} \R is free, it
suffices to conclude that for every $j < i$,
\begin{mathprop}{Quotient_j}%
        \Quotient {\S {\delta_{(j+1)}}} {(\S {\delta_j} \Union \R)} 
\end{mathprop}%
is free. Note that in this case $\delta_i$ might be in \LSE \emptyset
or in its complement, but in both cases \SfProj {\SimpleSeq
{\delta_j}} is a subset of \N i for every $j < i$
\Note {because of the inequality $\sigma > \lambda_1 > \dots >
\lambda_\rho$, $\Lh \rho > 1$}.
Since $(\delta_j) +1$ is not in \LSE \emptyset,
\Quotient {\S {\delta_{(j+1)}}} {\S {(\delta_j) +1}} is free.
Since \Quotient {\S {\delta_{(j+1)}}} {\S {(\delta_j) +1}} is in \N i,
also 
\[
        \Quotient
                {\Par[\big]{
                        \Quotient {\S {\delta_{(j+1)}}} {\S {(\delta_j) +1}}
                }}
                {\Par[\Big]{
                        \Par [\big] {
                                \Quotient {\S {\delta_{(j+1)}}} {\S {(\delta_j) +1}}
                        }
                        \Inter \N i
                }}
        =
        \Quotient
                {\S {\delta_{(j+1)}}}
                {\Par [\big] {
                        \S {(\delta_j) +1} \Union
                        (\R \Inter \S {(\delta_j) +1})
                }}
\]
is free.
Because both
$\BigUnion {
        (\S {(\delta_j) +1} \Inter \R)
}
= (\BigUnion {\S {(\delta_j) +1}})
\Inter
(\BigUnion \R)
$
and $\SFamilyProj {\SimpleSeq {\delta_j}} \Subset \R$ hold, we have
that
$
\Quotient {\S {\delta_{(j+1)}}} {(\S {\delta_j} \Union \R)} =
\Quotient {\S {\delta_{(j+1)}}} {
        \Par [\big] {\S {(\delta_j) +1}
        \Union (\R \Inter \S {(\delta_j) +1})}
}
$.
Hence \Property {Quotient_j} holds, and we have proved that \Quotient
{\S {\delta_i}} \R is free.


Let $X$ denote the set $\LSE \emptyset \Inter \delta_i \Minus \N
i$. Suppose, for the moment, that for every $\alpha \in X$, the set
\begin{textprop}{small}%
        $J_\alpha =
        \Set {\eta \in \SfProj {\SimpleSeq \alpha}}
        {\s \nu \Inter \s \eta \not= \emptyset
        \ForSome \nu \in \SfProj {\SimpleSeq {\delta_i}}}$
        is small in \SfProj {\SimpleSeq \alpha}.
\end{textprop}%
By \Property {Trans} choose a transversal $U_\alpha$ for $J_\alpha$,
$\alpha \in X$, so that $U_\alpha(\eta) \not\in \slevel \eta 1$ hold
for all $\eta \in J_\alpha$. By \RefOfDefinition {NPTSkeleton}
{det_by_index} together with the fact that $X \Inter \N i$ is empty,
$\Ran {U_\alpha} \Inter \BigUnion \R$ is empty. Consequently, the
desired extension $T'$ in \Property {below_delta_i} can be
$T \Union W \Union \BigUnion [\alpha \in X] {U_\alpha}$,
where $W$ is a restriction of some transversal witnessing that
\Quotient {\S {\delta_i}} \R is free into the set $\S {\delta_i}
\Minus \Set {\s \eta} {\eta \in J_\alpha \ForSome \alpha \in X}$.

So it remains to prove \Property {small}. Fix some \alpha from $X$.
Since $f_\alpha \EventLess f_{\delta_i}$, there is \Star \xi with
$f_\alpha(\xi) < f_{\delta_i}(\xi)$ for each $\xi > \Star \xi$.
Assume that $\eta \in \BigUnion [\alpha \in X] {\SfProj {\SimpleSeq
\alpha}}$ and $\nu \in \SfProj {\SimpleSeq {\delta_i}}$ are such that
$\s \eta \Inter \s \nu$ is nonempty. From \RefOfDefinition
{NPTSkeleton} {det_by_index} it follows that $\slevel \eta 1 \Inter
\slevel \nu 1 \not= \emptyset$ and there is \xi with $\eta(m) = \nu(m)
= \xi$. From the new definition of the ``first coordinate'' it follows
that $f_\alpha(\xi) = f_{\delta_i}(\xi)$. Therefore $\xi \leq \Star
\xi$. This means that for every $\tau \in \LS$ of length $m$, the set
\Set {\eta \in \SfProj \tau} {\s \eta \Inter \s \nu \not= \emptyset
\ForSome \nu \in \SfProj {\SimpleSeq {\delta_i}}} has cardinality $<
\lambda_{\Res \eta m} = \lambda_{\Res \nu m} = \Cf \mu$. Thus
\Property {small} holds and we have proved \Proposition {Singular}.

%

\end{SECTION}


\begin{SECTION}{NRT}{Building blocks of incompactness skeletons}

%
In this section $\lambda > \kappa$ are regular cardinals and \seq
\kappa is a decreasing sequence of regular cardinals.
In the definitions below we present a variant ``\NRT \lambda \kappa''
of the notion \NPT \lambda \kappa. A simplest form of this variant is
mentioned in \cite [Sh355: Fact 6.2(9)] {CardArith}.
So we define analogical notions of ``a transversal for a subfamily''
and ``a free subfamily over another subfamily''.
There are two main motivations. Firstly, even if both \NPT {\kappa_1}
{\kappa_2} and \NPT {\kappa_2} {\kappa_3} hold, for regular cardinals
$\kappa_1 > \kappa_2 > \kappa_3$, \NPT {\kappa_1} {\kappa_3} does not
need to hold:

\begin{FACT}{NotTrans}%
\AddInTheoremCite[Lemma 3.1]{Sh523}%
It is consistent with \ZFC relative to existence of two weakly compact
cardinals, that for some regular cardinals $\lambda > \kappa > \Alephf
0$, both \NPT \lambda \kappa and \NPT \kappa - hold, but \NPT \lambda
- does not hold.
\end{FACT}

However, we show that an analogical ``transitivity property'' for \NRT
\lambda \kappa holds, \Corollary {Transitive}.
Secondly, the introduced families and their canonical forms, called
\NRT \lambda {-}-skeletons, provide a unified picture of those
construction methods of \NPT \lambda {-}-skeletons presented in \cite
[\SS1] {MgSh204} and \cite [\SS2] {Sh523} \Note {see \Fact
{ExistenceOfNRT}}.
We shall deal with the ideal \NonstatIdeal {\seq \kappa} of ``the
product'' of nonstationary ideals of fixed cardinals \seq \kappa \Note
{the next definition} instead of considering just the ideal of the
bounded subsets of \kappa.
However, the reader may think in the beginning, if she or he wants,
that $\seq \kappa = \SimpleSeq \kappa$ and \seq \Ideal consist of the
ideal of all bounded subsets of \kappa.

\begin{REMARK}
\relax From \Fact {NPTImpliesNRT} below it follows that \NRT \lambda {\seq
  \kappa}-skeletons are necessary to ``pump up'' more complicated \NPT
\lambda {-}-skeletons.
In fact we shall see in \Definition {NRTSet} and \Conclusion
{NicelyIncompact}, that the definition of a nicely incompact set of
regular cardinals from \cite [Definition 2.1] {Sh521} can be expressed
using \NRT \lambda {\seq \kappa}-skeletons as building blocks.
All together \Note {which means \Fact {AddIndex}, \Definition
{Compatible}, \Lemma {Transitive}, and \Lemma {AddLevel}}, \NRT
\lambda {\seq \kappa}-skeletons are handy tools to divide and
amalgamate \NPT \lambda {-}-skeletons.
The reader is advised to look at \Section {Conclusion} to understand
our goal, if she or he feel lost with the technicalities in this
section.
\end{REMARK}

\begin{DEFINITION}{AFamily}
Suppose $\seq \kappa = \Seq {\kappa_k} {k < \Star k}$ is a nonempty
decreasing sequence of infinite cardinals.
A family \AFamily is called an \NRT-{\seq \kappa}-family
\Note {\NRT-{\seq \kappa} comes from the notion \NRT \lambda {\seq
\kappa}},
when \AFamily fulfills the following demands:
\begin{properties}

\item the elements of \AFamily are finite sequences of functions;

\item \AFamily is enumerated by \Set {\ASeq i} {i < \Card \AFamily};

\item each \ASeq i is of the form \Seq {a_{i,l}} {l < \Lh {\ASeq i}};

\item $a_{i,l}$'s have domain \Prod {\seq \kappa}.

\end{properties}
\end{DEFINITION}

\begin{DEFINITION}{ProductOfIdeals}
Suppose \seq \kappa is as in \Definition {AFamily}.
A sequence $\seq \Ideal = \Seq {\Ideal_k} {k < \Star k}$ is called a
\seq \kappa-sequence of ideals, when each $\Ideal_k$ is a proper ideal
on $\kappa_k$ containing all the bounded subsets of $\kappa_k$.
For two ideals $\Ideal_0$ and $\Ideal_1$ on $\kappa_0$ and $\kappa_1$
respectively, the product $\Ideal_0 \times \Ideal_1$ is the ideal of all
subsets $X$ of $\kappa_0 \times \kappa_1$ such that
\[
        \Set [\Big] {\alpha_0 < \kappa_0}
        {
                \Set [\big] {\alpha_1 < \kappa_1} {
                        \Pair {\alpha_0}{\alpha_1} \in X
                } \not\in \Ideal_1
        } \in \Ideal_0.
\]
For a regular cardinal \kappa, let \BoundedIdeal \kappa denote the
ideal of all the bounded subsets of \kappa, and if \kappa is
uncountable, let \NonstatIdeal \kappa denote the ideal of all
nonstationary subsets of \kappa.
For a \seq \kappa-sequence of ideals, let \NonstatIdeal {\seq \kappa}
denote the product of the ideals in \Seq {\Ideal_k'} {k < \Star k},
where $\Ideal_k' = \NonstatIdeal {\kappa_k}$ if $\kappa_k$ is
uncountable, and otherwise, $\Ideal_k = \BoundedIdeal {\Alephf 0}$.
\end{DEFINITION}

\begin{DEFINITION}{ATrans}
Assume that \seq \Ideal is a \seq \kappa-sequence of ideals, \AFamily
is an \NRT-{\seq \kappa}-family, and $\AFamily['] = \Set {\ASeq i} {i
  \in I}$ is a subfamily of \AFamily.
A sequence $\ATrans b = \Seq {b_i} {i \in I}$ of functions is called
an \NRT[\seq \Ideal]-{\seq \kappa}-transversal for \AFamily['] if the
following properties are satisfied for every $i \in I$:
\begin{properties}

\item for some $l < \Lh {\ASeq i}$, $b_i \Subset a_{i,l}$ and $\Prod
{\seq \kappa} \Minus \Dom {b_i}$ is in \Prod {\seq \Ideal};

\item the ranges of $b_i$, $i \in I$, are pairwise disjoint.

\end{properties}
\AFamily['] is called \NRT[\seq \Ideal]-{\seq \kappa}-free if there
exists an \NRT[\seq \Ideal]-{\seq \kappa}-transversal for \AFamily[']
\Note {the empty sequence is transversal for the empty family}.
A family \AFamily is almost \NRT[\seq \Ideal]-{\seq \kappa}-free, when
every subfamily of cardinality $< \Card {\AFamily}$ is \NRT[\seq
\Ideal]-{\seq \kappa}-free.

For a cardinal $\lambda > \max \seq \kappa$, we say that \NRT [\seq
\Ideal] \lambda {\seq \kappa} holds, if there exists an almost \NRT
[\seq \Ideal]-{\seq \kappa}-free \NRT-{\seq \kappa}-family of
cardinality \lambda which is not \NRT [\seq \Ideal]-{\seq
\kappa}-free.
In the case $\seq \kappa = \SimpleSeq \kappa$ and $\seq \Ideal =
\SimpleSeq {\BoundedIdeal \kappa}$, \NRT \lambda \kappa is an
abbreviation for \NRT [\seq \Ideal] \lambda {\seq \kappa}.
\end{DEFINITION}

First we see that existence of an \NPT \lambda {-}-skeleton
necessarily gives many almost free nonfree \NRT-{\seq \kappa}-families
\Note {depending on the type of the skeleton}.
Later in \Lemma {Transitive} we shall see that existence of certain
type of \NRT-{\seq \kappa}-families is also sufficient condition for
building \NPT \lambda {-}-skeletons.
More generally, ``neat'' \NRT-{\seq \kappa}-families can be
transformed into a form of an ``\NRT \lambda {\seq \kappa}-skeleton'',
which is an analog of \NPT \lambda {-}-skeletons for \NRT-{\seq
  \kappa}-families.

\begin{FACT}{NPTImpliesNRT}%
Let \NPTSkeleton be an \NPT \lambda {-}-skeleton of type \seq \lambda
and of height \Star n.
Suppose that there is a sequence \Seq {\lambda_n} {n < \Star n} such
that for every $n < \Star n$ and $\rho \in \LS$ of length $n$,
$\lambda_n = \lambda_\rho$ holds \Note {\seq \lambda does not contain
regular limit cardinals}.
Then for every $n < \Star n$, \NRT [\NonstatIdeal {{\seq \kappa}^n}]
{\lambda_n} {{\seq \kappa}^n} holds, where ${\seq \kappa}^n =
\ConcatSimple {\Seq {\lambda_m} {n < m < \Star n}} {\Alephf 0}$ and
\NonstatIdeal {{\seq \kappa}^n} is the ideal given in \Definition
{ProductOfIdeals} above.
In fact, there exists an \NRT {\lambda_n} {\seq \kappa}-skeleton of
type \SimpleSeq {\lambda_n}, see \Definition {NRTSkeleton}.

\begin{Proof}%
Fix $n < \Star n$ and $\rho \in \LS \Minus \Sf$ of length $n$.
Let \Star k denote the length of \seq \kappa, i.e., $\Star k = \Star n
-n +1$.
Define a family \AFamily [_\rho] to be \Set [\big] {\SimpleSeq
{a_\alpha}} {\alpha \in \LSE \rho}, where each $a_\alpha$ is a
function having domain \Prod {\seq \kappa} and for every $\tau \in
\Prod {\seq \kappa}$, the value $a_\alpha(\tau)$ is chosen in the
following way.
For fixed \alpha and \tau let $\eta_{\alpha, \tau} \in \Sf$ be such
that $\Res {\eta_{\alpha, \tau}} n = \rho$, $\eta_{\alpha, \tau}(n) =
\alpha$, and for every $k < \Star k - 1$, $\eta_{\alpha, \tau}(n+k)$
is the $\tau(k)$'s member of \LSE {\Res {\eta_{\alpha, \tau}} {n+k}}
in $\in$-order.
Define $a_\alpha(\tau)$ to be the $\tau(\Star k-1)$'s element of the
countable set $\slevel {\eta_{\alpha, \tau}} {n+1} = \s {\eta_{\alpha,
\tau}} \Inter \LSet {\ConcatSimple \rho \alpha}$ in the fixed
enumeration of it, recall \RefOfDefinition {NPTSkeleton} {tree-order}.

By Fodor's lemma \AFamily [_\rho] cannot be \NRT[\NonstatIdeal {\seq
\kappa}]-{\seq \kappa}-free. On the other hand \AFamily [_\rho] is
almost \NRT[\NonstatIdeal {\seq \kappa}]-{\seq \kappa}-free:
Fix $\delta < \lambda_n$ and let $I$ denote the set \Set {\eta \in
\SfProj \rho} {\eta(n) < \delta}. By \cite [Claim 3.8(D)] {Sh161} we
may choose a sequence \Seq {u_\eta} {\eta \in I} of pairwise disjoint
sets such that for every $\eta \in I$ there is an index $l_\eta$ such
that $n \leq l_\eta < \Star n$ and $u_\eta$ is an end segment of
\slevel \eta {l_\eta + 1}. For every $\alpha < \delta$ the set
\[
        I^\alpha = \Set [\big] {\nu \in \SfProj \rho}
        {\nu(n) = \alpha \And l_\nu > n}
\]
must be small in \SfProj {\ConcatSimple \rho \alpha}.
It follows that the corresponding set $Y^\alpha = \Set {\tau \in \Prod
{\seq \kappa}} {\eta_{\alpha, \tau} \in I^\alpha}$ must be in
\NonstatIdeal {\seq \kappa}. So we may define a transversal \Seq
{b_\alpha} {\alpha < \delta} for \Set [\big] {\SimpleSeq {a_\alpha}}
{\alpha < \delta} by setting $b_\alpha = \Res {a_\alpha} {\Prod {\seq
\kappa} \Minus Y^\alpha}$.
\end{Proof}%
\end{FACT}

There are many examples of nonfree almost free \NRT [\seq
\Ideal] {-} {\seq \kappa}-families.
Note that the examples below are \NRT \lambda {\seq \kappa}-skeletons
of height 1 for some \seq \kappa of length at most 2.

\begin{FACT}{ExistenceOfNRT}%
Suppose \lambda is an uncountable regular cardinal.
\begin{ITEMS}

\ITEM{Nonreflecting}%
\AddInTheoremCite[Lemma 23]{Sh108}%
%
%
%
If there exist a regular cardinal $\kappa < \lambda$ and a
nonreflecting stationary subset $E$ of $\lambda \Inter \Cof \kappa$
\Note {which means for all $\alpha < \lambda$, $E \Inter \alpha$ is
nonstationary in \alpha},
then \NRT \lambda \kappa holds.

\ITEM{Successor}%
If $\lambda = \SuccCard \kappa$ and \kappa is a regular cardinal, then
\NRT \lambda \kappa holds \Note {since $\lambda \Inter \Cof \kappa$ is
a nonreflecting stationary set, or for the other cofinalities use
$*(\kappa, \kappa)$ given in the references of the next item}.

\ITEM{StarPrinc}%
\AddInTheorem{%
\cite{MkSh251} or %
\cite[Theorem VI.3.9 and VI.3.10]{EklofMekler}%
}
If \mu is a singular cardinal of cofinality \kappa and \RealIsConstr,
then \NRT {\SuccCard \mu} \kappa holds \Note {because of $*(\mu,
\kappa)$}.

\ITEM{Sh108}%
\AddInTheoremCite[Lemma 24]{Sh108}%
Suppose \mu is a singular strong limit cardinal of cofinality \kappa
such that $\IGoodIdeal {\SuccCard \mu} = \PowerSet {\SuccCard \mu}$
\Note {for a definition of \IGoodIdeal {\SuccCard \mu}, see e.g. \cite
[Definition 2.1] {Sh351}}.
Then for every regular $\theta < \mu$ with $\kappa \not= \theta$,
\NRT [\seq \Ideal] \lambda {\Pair \theta \kappa} holds,
where $\seq \Ideal = \Pair {\BoundedIdeal \theta} {\BoundedIdeal
\kappa}$.

\ITEM{MgSh204}%
\AddInTheoremCite[Theorem 4]{MgSh204}%
Suppose $\kappa < \theta$ are regular cardinals and $\lambda =
\SuccCard [\kappa+1] \theta$.
Then \NRT [\seq \Ideal] \lambda {\Pair \theta \kappa} holds,
where $\seq \Ideal = \Pair {\NonstatIdeal \theta} {\BoundedIdeal
\kappa}$.

\ITEM{Sh523A}%
\AddInTheoremCite[Lemma 1.16]{Sh523}%
Suppose \kappa is a regular cardinal such that for $\mu = \SuccCard
[\kappa] \kappa$ there is a scale \Pair {\seq \mu} {\scale f} whose
good points \good {\scale f} contains a closed and unbounded subset of
\SuccCard \mu \Note {\Definition {scale} and \Definition {good}}.
Then for every regular $\theta$ with $\kappa < \theta < \mu$, \NRT
[\seq \Ideal] {\SuccCard \mu} {\Pair \theta \kappa} holds,
where $\seq \Ideal = \Pair {\NonstatIdeal \theta} {\BoundedIdeal
\kappa}$.

%
%
%

\ITEM{Sh355}%
\AddInTheoremCite[Sh355:Claim II.1.5A]{CardArith}%
If \mu is a singular cardinal of cofinality \kappa and $\pp
[^*_\Ideal] \mu > \SuccCard \mu$,
then \NRT [\Ideal] {\SuccCard \mu} \kappa holds.
Moreover, when \kappa is uncountable, already $\pp \mu > \SuccCard
\mu$ together with some weak assumptions imply \NRT [\BoundedIdeal
\kappa] {\SuccCard \mu} \kappa
\cite [Analytical Guide 5.7(B) and Sh371:\SS0 and \SS1] {CardArith}.

\ITEM{Sh355B}%
\AddInTheoremCite[Sh355:Theorem II.6.3]{CardArith}%
If \mu is a singular cardinal of countable cofinality and
$\Symb{cov}(\mu, \mu, \Alephf 1, 2) > \SuccCard \mu$ \Note {e.g., \mu
strong limit and $\Exp \mu {\Alephf 0} > \SuccCard \mu$}, then \NRT
{\SuccCard \mu} {-} holds.

\end{ITEMS}

\end{FACT}

The rest of this section is essentially based on a similar analysis
and a construction of a \lambda-system as in \cite [\SS3] {Sh161}.
Our presentation resembles \cite [Appendix] {Sh161} \Note {or \cite
[VII.3A] {EklofMekler}}.
The next definition yields an analog of the notion ``a subfamily is
free over an other subfamily''.

\begin{DEFINITION}{Quotient}%
Suppose \AFamily is an \NRT-{\seq \kappa}-family.
A sequence \Seq {\A \alpha} {\alpha < \lambda} is called a filtration
of \AFamily if it is a continuous increasing chain of subfamilies,
such that the members have cardinality $< \lambda$ and the union of
the members equals \AFamily.

For every $\AFamily['] \Subset \AFamily$, \AUnion {\AFamily[']}
denotes the set \BigUnion [\seq a \in \AFamily'] {\BigUnion [l < \Lh
{\seq a}] {\Ran {a_l}}}.

Suppose \seq \Ideal is a fixed \seq \kappa-sequence of ideals and \A
1, \A 2 are subfamilies of \AFamily. 
We denote by \Quotient {\A 2} {\A 1} the family \Set {\Res {\seq a}
{L_{\seq a}}} {\seq a \in \A 2 \Minus \A 1 \And L_{\seq a} \not=
\emptyset}
\Note {\seq \Ideal should be clear from the context}, where $L_{\seq
a}$ denotes those indices $l < \Lh {\seq a}$ for for which $\Ran {a_l}
\Inter (\AUnion {\A 1})$ is ``small'', i.e.,
\[
        L_{\seq a} =
        \Set [\Big] {l < \Lh {\seq a}} {
                \Set [\big] {\tau \in \Prod {\seq \kappa}}
                {a_l(\tau) \in \AUnion {\A 1}} \in \Prod {\seq \Ideal}
        }.
\]
We say that \A 2 is free over \A 1, if the family \Quotient {\A 2} {\A
1} is free \Note {we omit the prefix \NRT[\seq \Ideal]-{\seq \kappa}}.
\end{DEFINITION}

The next fact offers some basic facts needed to understand \Lemma
{SingComp}.

\begin{FACT}{Cub}%
Suppose \AFamily is an \NRT-{\seq \kappa}-family and \seq \Ideal is a
\seq \kappa-sequence of ideals.
\begin{ITEMS}

\ITEM{free}%
If \Seq {\A \alpha} {\alpha < \delta} is a continuous chain of
subfamilies of \AFamily such that \A 0 is free and \Quotient {\A
{\alpha+1}} {\A \alpha} is free for every $\alpha < \delta$, then the
family \BigUnion [\alpha < \delta] {\A \alpha} is free.

\ITEM{regular}%
Suppose \AFamily is of regular cardinality \lambda and \AFamily is
free. Then for all filtrations \Seq {\A \alpha} {\alpha < \lambda} of
\AFamily, the set \Set {\alpha < \lambda} {\Quotient {\AFamily} {\A
\alpha} \Text{is free}} contains a cub of \lambda.

\end{ITEMS}

\begin{Proof}%
\ProofOfItem{free}%
Build an \NRT[\seq \Ideal]-{\seq \kappa}-transversal for the union by
induction on $\alpha < \delta$ using the following: Suppose $\A 1
\Subset \A 2 \Subset \A 3$ are subfamilies of \AFamily, \ATrans b is
an \NRT[\seq \Ideal]-{\seq \kappa}-transversal for \Quotient {\A 2}
{\A 1}, and \ATrans c is an \NRT [\seq \Ideal]-{\seq
  \kappa}-transversal for \Quotient {\A 3} {\A 2}. Then the
concatenation \Concat {\ATrans b} {\ATrans c} is an \NRT [\seq
\Ideal]-{\seq \kappa}-transversal for \Quotient {\A 3} {\A 1}, since
w.l.o.g. for every $i \in I$, the intersection $\Ran {b_i} \Inter
(\AUnion {\A 1})$ is empty.

\ProofOfItem{regular}%
Suppose \Seq {b_i} {i < \lambda} is an \NRT[\seq \Ideal]-{\seq
  \kappa}-transversal for \AFamily. By the demands on \seq \kappa and
\seq \kappa-sequence of ideals, each of the sets \Ran {b_i} has
cardinality $\kappa_0 = \max \seq \kappa$.
Let $h_\xi$, $\xi < \kappa_0$, be injective functions with domain
\lambda such that for every $i < \lambda$, $h_\xi(i)$ is the \xi's
element of \Ran {b_i} in some fixed enumeration.
To prove the claim, it suffices to choose a cub \C of \lambda so that
for all $\alpha \in \C$ and $i < \lambda$:
\begin{property}

if $h_\xi(i) \in \AUnion {\A \alpha}$ for some $\xi < \kappa_0$, then
$\ASeq i \in \A \alpha$.

\end{property}
Then \Seq {b_i} {i < \lambda \And \ASeq i \not\in \A \alpha} is an
\NRT[\seq \Ideal]-{\seq \kappa}-transversal for \Quotient {\AFamily}
{\A \alpha}.
\end{Proof}%
\end{FACT}

\begin{LEMMA}{SingComp}%
Suppose \seq \kappa and \seq \Ideal are as in \Definition
{ATrans}. For all singular cardinals $\mu > \max \seq \kappa$, every
almost \NRT[\seq \Ideal]-{\seq \kappa}-free family of cardinality \mu
is free, i.e., \NRT [\seq \Ideal] \mu {\seq \kappa} does not hold.

\begin{Proof}%
The claim follows from \cite [Theorem 2.1] {Sh52} or \cite [Theorem in
\SS0] {BenDavid1978}, because for a fixed family \AFamily of
cardinality \lambda, the relation ``\A 2 is free over \A 1'' for
subfamilies of \AFamily satisfies the demanded axioms. However, we
briefly sketch why axioms I--V of \cite [Theorem 5 in \SS4]
{Hodges1981} hold.
Define the set $S(\AFamily)$ of ``the subalgebras of \AFamily'' to be
the set of all subfamilies of \AFamily. A subalgebra $\A 1 \in
S(\AFamily)$ is free when \A 1 has an \NRT [\seq \Ideal]-{\seq
  \kappa}-transversal, say \ATrans b, and a basis $\mathcal F$ of \A 1
is the set of all subalgebras \A 2 of \A 1 such that the corresponding
restriction of \ATrans b is a transversal for \Quotient {\A 1} {\A 2}
\Note {the proof of \ItemOfFact {Cub} {regular}}. Note that a basis
$\mathcal F$ is ``fully closed unbounded above $\kappa_0$'' as
demanded in axiom II \Note {use $h_\xi$'s as in the proof of
  \ItemOfFact {Cub} {regular}}. Axioms I, III, and IV hold by the
definition. Axiom V can be proved by a similar construction as in the
proof of \ItemOfFact {Cub} {free}.
\end{Proof}
\end{LEMMA}

Our next task is to prove that \NRT {-} {\seq \kappa}-families can be
transformed into a canonical form in the same way as families of
countable sets are transformed into incompactness skeletons \Note
{\Lemma {NRTSkeleton}}. In order to succeed in the proof we have to
assume the \NRT {-} {\seq \kappa}-family to be ``neat'' \Note
{\Definition {NeatAFamily}}.

\begin{DEFINITION}{NRTSkeleton}%
Suppose \seq \kappa and \seq \Ideal are as in \Definition {ATrans} and
\lambda is a regular cardinal greater than $\max \seq \kappa$.
A tuple \NRTSkeleton is called an \NRT \lambda {\seq \kappa}-skeleton
of type \seq \lambda, when the following conditions hold
\Note {recall the definitions of a canonical \lambda-set and a
disjoint \lambda-system from \Section {Preliminaries}}:
\begin{properties}

\item \LS is a \lambda-set of type $\seq \lambda = \Seq {\lambda_\rho}
{\rho \in \LS \Minus \Sf}$;

\item \LS has a canonical form and its height is $\Star n < \omega$;

\item $\LSystem = \Seq [\big] {\LSet {\ConcatSimple \rho \alpha}}
{\rho \in \LS \Minus \Sf \And \alpha < \lambda_\rho}$ is a disjoint
\lambda-system;

\item \AFamily is an \NRT-{\seq \kappa}-family of cardinality \lambda;

\item \AFamily is enumerated by \Set {\ASeq \eta} {\eta \in \Sf};

\item every sequence in \AFamily has a fixed length $\Star l <
\omega$, where $\Star l \geq \Star n$;

\item \AFamily is based on \LSystem, i.e., for all $\eta \in \Sf$ and
$l < \Star l$, $\Ran {a_{\eta, l}} \Subset \BigUnion [m \leq \Lh \eta]
{\LSet {\Res \eta m}}$;

\item \AFamily is almost \NRT [\NonstatIdeal {\seq \kappa}]-{\seq
\kappa}-free
\Note {\NonstatIdeal {\seq \kappa} given in \Definition
{ProductOfIdeals}}.

\end{properties}
Moreover, analogously to the definition of an \NPT \lambda
\kappa-skeleton, we demand that for every $n < \Star n$ and $\rho \in
\LS$ of length $n$ the following conditions hold:

When the cardinal $\theta_n$, given by the canonical form of \LS \Note
{\Definition {CanonicalLambdaSet}}, is well-defined:
\begin{properties}

\item $\theta_n > \max \seq \kappa$ implies that there is $m < \Star
n$ such that the cardinal $\lambda_m$, given by the canonical form of
\LS, is well-defined and $\lambda_m = \theta_n$;

\item if $\Alephf 0 < \theta_n \leq \max \seq \kappa$ then $\theta_n
\in \Ran {\seq \kappa}$.

\end{properties}

The index set \Star l can be partitioned into \Star n blocks \Seq
{L^{n+1}} {n < \Star n}
\Note {analogously to the partition \Seq {\slevel \eta {n+1}} {n <
    \Star n} of a set in an \NPT \lambda {-}-skeleton}
such that for every $n < \Star n$, $l \in L^{n+1}$, and $\eta \in
\Sf$,
        $\Ran {a_{\eta,l}} \Subset \LSet {\Res \eta {n+1}}$
\Note {by the disjointness of \LSystem, for all $l \in L^{n+1}$ and
$l' \in L^{n'+1}$, $\Ran {a_{\eta,l}} \Inter \Ran {a_{\eta',l'}} =
\emptyset$ whenever $\Res \eta n \not = \Res {\eta'} n$}.

For all \rho in $\LS \Minus \Sf$, \NRTSkeletonProj \rho denotes the
\NRT {\lambda_\rho} {\seq \kappa}-skeleton of type \Seq {\lambda_\tau}
{\tau \in \LSProj \rho \Minus \Sf}, where
\begin{properties}

\item \LSProj \rho denotes the set \Set {\eta \in \LS} {\rho \Subset
\eta};

\item \LSystemProj \rho is the restriction of \LSystem according to
the new index set;

\item \AFamilyProj \rho denotes the family \Set [\big] {\Res {\ASeq
\eta} {L^{>\Lh \rho}}} {\eta \in \SfProj \rho}, where
\[
L^{>\Lh \rho} = \BigUnion [\Lh \rho \leq n < \Star n] {L^{n+1}}.
\]

\end{properties}

Lastly, for every $\rho \in \LS$ of length $n > 0$ and $\gamma <
\rho(n-1)$, there are strictly less than $\lambda_\rho$ indices $\eta
\in \Sf$ such that $\rho \Subset \eta$ and for some $l \in L^{n+1}$,
the set \Set [\big] {\tau \in \Prod {\seq \kappa}} {a_{\eta,l}(\tau)
\in \LSet {\ConcatSimple \rho \gamma}} is not in \NonstatIdeal {\seq
\kappa}.
\end{DEFINITION}

\begin{FACT}{NotFree}%
If \NRTSkeleton is an \NRT \lambda {\seq \kappa}-skeleton, then there
is no injective choice function for the family \Set {\s \eta} {\eta
  \in \Sf}, where each \s \eta is the set \BigUnion [l < \Lh {\ASeq
  \eta}] {\Ran {a_{\eta,l}}}.
In particular, there is no \NRT [\NonstatIdeal {\seq \kappa}]-{\seq
\kappa}-transversal for \AFamily.

\begin{Proof}%
The given family is based on the \lambda-system \LSystem. Hence the
claim follows from \cite [Claim 3.5] {Sh161} \Note {or \cite [VII.3.6]
{EklofMekler}}.
\end{Proof}%
\end{FACT}

\begin{DEFINITION}{NeatAFamily}%
For technical reasons \Note {needed in the proof of \Lemma
  {NRTSkeleton}}, we say that an \NRT-{\seq \kappa}-family \AFamily is
neat, when for all $\seq a, \seq a' \in \AFamily$, $l < \Lh {\seq a}$,
and $l' < \Lh {\seq a'}$:
$a_l \not= a_{l'}'$ implies that the set
\[
        \Set [\big] {\tau(0)} {\tau \in \Prod {\seq \kappa}
        \And a_l(\tau) = a_{l'}'(\tau)}
\]
is a nonstationary subset of $\kappa_0$, if $\kappa_0$ is uncountable,
and otherwise, the set is finite.
\end{DEFINITION}

\begin{REMARK}
We do not demand neatness in the definition of an \NRT \lambda {\seq
  \kappa}-skeleton. Therefore \ItemOfFact {ExistenceOfNRT} {Sh108},
\ItemOfFact {ExistenceOfNRT} {MgSh204}, and \ItemOfFact
{ExistenceOfNRT} {Sh523A} are examples of \NRT \lambda {\seq
  \kappa}-skeletons.
\end{REMARK}

%
%
In the definition of the neatness the most natural demand would be
that ``intersection of two different coordinates'' is small in the
\NonstatIdeal {\seq \kappa} sense, but such a definition would cause
problems in the next lemma, since \NonstatIdeal {\seq \kappa} is not
$\kappa_0$-complete when \seq \kappa has length greater than 1.

\begin{LEMMA}{NRTSkeleton}%
If \AFamily is a neat family witnessing that \NRT [\NonstatIdeal {\seq
  \kappa}] \lambda {\seq \kappa} holds, then there is an \NRT \lambda
{\seq \kappa}-skeleton.

\begin{Proof}
%
We define an \NRT \lambda {\seq \kappa}-skeleton \NRTSkeleton[']
following the same ideas as in \cite [Claim 3.3] {Sh161} or \cite
[Appendix: Proposition 4] {Sh161} \Note {or \cite [Proposition
  VII.3.7] {EklofMekler}}, except that at the end we need a different
trick.

Let $\lambda_\emptyset$ be \lambda itself.  By \ItemOfFact {Cub}
{free} there is a filtration \Seq {\A {\SimpleSeq \alpha}} {\alpha <
  \lambda_\emptyset} of \AFamily and a stationary subset
$\E_\emptyset$ of \lambda such that for every $\alpha \in
\E_\emptyset$, \Quotient {\A {\SimpleSeq {\alpha+1}}} {\A {\SimpleSeq
    \alpha}} is not free. For every $\alpha \in \E_\emptyset$, choose
$\B {\SimpleSeq \alpha} \Subset \A {\alpha+1} \Minus \A \alpha$ so
that $\lambda_{\SimpleSeq \alpha} = \Card {\B {\SimpleSeq \alpha}} <
\lambda_\emptyset$ is the smallest cardinal such that \B {\SimpleSeq
  \alpha} is not free over \A \alpha.

Suppose $\lambda_{\SimpleSeq \alpha} > \kappa_0$ holds. By the choice
of \B {\SimpleSeq \alpha}, the family \Quotient {\B {\SimpleSeq
\alpha}} {\A {\SimpleSeq \alpha}} witnesses that \NRT [\NonstatIdeal
{\seq \kappa}] {\lambda_{\SimpleSeq \alpha}} {\seq \kappa} holds. By
\Lemma {SingComp} $\lambda_{\SimpleSeq \alpha}$ is a regular cardinal.
By \ItemOfFact {Cub} {free} there is a filtration \Seq {\A {\SimpleSeq
{\alpha, \beta}}} {\beta < \lambda_\alpha} of \B {\SimpleSeq \alpha}
and a stationary subset $\E_{\SimpleSeq \alpha}$ so that
\A {\SimpleSeq {\alpha, \beta+1}} is not free over $\A {\SimpleSeq
{\alpha, \beta}} \Union \A {\SimpleSeq {\alpha}}$. Hence we may
continue choosing a subset \B {\SimpleSeq {\alpha, \beta}} of $\A
{\SimpleSeq {\alpha, \beta+1}} \Minus (\A \alpha \Union \A {\SimpleSeq
{\alpha, \beta}})$ having the smallest possible cardinality
$\lambda_{\SimpleSeq {\alpha,\beta}} < \lambda_{\SimpleSeq \alpha}$
for which \B {\SimpleSeq {\alpha, \beta}} is not free over $\A
{\SimpleSeq {\alpha, \beta}} \Union \A {\SimpleSeq {\alpha}}$.

Assume \eta is a sequence of ordinals such that its length is $n > 0$
and $\lambda_{\Res \eta {n-1}} > \kappa_0 \geq \lambda_\eta$.  Such
sequences \eta form the final nodes \Sf['] of the desired \lambda-set
\LS[']. Suppose also that all the sets \A {\Res \eta m}, for nonzero
$m \leq n$, together with \B \eta are already chosen so that \B \eta
is not free over \A [^*] \eta, where \A [^*] \eta is an abbreviation
for \BigUnion [0 < m \leq n] {\A {\Res \eta m}}. To simplify our
explanations, let \Quotient {\seq a} {\A [^*] \eta} denote the
sequence \Note {possibly empty} in the singleton \Quotient {\Braces
{\seq a}} {\A [^*] \eta}, and on the other hand, let $\seq a \AInter
\A [^*] \eta$ denote the ``complementary'' sequence \Seq [\big] {a_l}
{l < \Lh {\seq a} \And l \not\in \Dom {\Quotient {\seq a} {\A [^*]
\eta}}}.
Now we would like to define that the new \NRT-{\seq \kappa}-family is
$\AFamily['] = \Set {\ASeq \eta'} {\eta \in \Sf[']}$, where each
$\ASeq \eta'$ is $\seq a \AInter \A [^*] \eta$ for some $\seq a \in \B
\eta$. However, we should choose $\ASeq \eta'$ so that \AFamily[']
becomes based on a \lambda-system.

The desired \lambda-system \LSystem['] is defined by setting for every
$\rho \in \LS['] \Minus \Sf[']$ and $\alpha < \lambda_\rho$, $\LSet[']
{\ConcatSimple \rho \alpha} = \AUnion {\A {\ConcatSimple \rho
\alpha}}$.

Assume that \AFamily is enumerated by \Set {\ASeq i} {i < \lambda}.
Let for every $i < \lambda$, $X_{\eta,i,l}$ be those indices for which
$a_{i,l}$ takes value in the defined \lambda-system, i.e.,
\[
        X_{\eta,i,l} =
        \Set [\big] {\tau \in \Prod {\seq \kappa}}{
                a_{i,l}(\tau) \in \AUnion {\A [^*] \eta}
        }.
\]
We get a better candidate for $\ASeq \eta'$ by considering the
sequence \Seq [\big] {\Res {a_{i,l}} {X_{\eta,i,l}}} {l < \Lh {\seq a}
\And l \not\in \Dom {\Quotient {\seq a} {\A [^*] \eta}}}.
Even this sequence is problematic. The domains of the coordinates are
not equal to \Prod {\seq \kappa} as demanded in the definition of an
\NRT-{\seq \kappa}-family.
\Note {Why the domain should be of that form? Otherwise we run into
difficulties when we combine two skeletons together in \Lemma
{Transitive}.}
We can fix this in the same way as we did in the proof of \Fact
{NPTImpliesNRT}. Let $\pi_{\eta,i,l}$ be the following ``continuous''
map from \Prod {\seq \kappa} onto $X_{\eta,i,l}$ defined for every
$\sigma \in \Prod {\seq \kappa}$ by
\begin{property}

        $\pi_{\eta,i,l}(\sigma) = \tau$ iff for each $k < \Lh {\seq
        \kappa}$, $\tau(k)$ is the $\sigma(k)$'s ordinal \Note {in
        $\in$-order} of the set \Set {\tau'(k)} {\tau' \in
        X_{\eta,i,l} \And \EqualRes {\tau'} \tau k}.
        
\end{property}
The property of these maps used below is that if $X_{\eta,i,l}$ is not
in \NonstatIdeal {\seq \kappa}, then for every $B \Subset \Prod {\seq
\kappa}$:
\begin{mathprop}{PiMap}
        \Prod {\seq \kappa} \Minus B \in \NonstatIdeal {\seq \kappa}
        \Iff \InvImage {\pi_{\eta,i,l}} {X_{\eta,i,l} \Minus B} \in
        \NonstatIdeal {\seq \kappa}.
\end{mathprop}%
A new candidate for $\ASeq \eta'$ could be \Seq {\Res {a_{i,l}}
  {X_{\eta,i,l}} \Comp \pi_{\eta,i,l}} {a_{i,l} \in \ASeq i \AInter \A
  [^*] \eta} for some $i$ such that $\ASeq i \in \B \eta$.  But then
we face ``the real problem'', why should such an \AFamily['] be almost
free?  A final definition of $\ASeq \eta'$ will be a concatenation of
finite number of sequences of the lastly given form. How to choose a
suitable finite set?

We claim that the neatness of the original \AFamily guarantees that
there is no ``choice function'' picking different functions from the
sequences in \Quotient {\B \eta} {\A [^*] \eta}, i.e., there is no
injective function $f$ with domain \Quotient {\B \eta} {\A [^*] \eta}
such that $f(\seq a) \in \seq a$ for each \seq a in the domain.
Namely, assume that such a $f$ exists and \Seq {\ASeq \xi} {\xi <
\lambda_\eta} enumerates \B \eta. By induction on $\xi < \lambda_\eta$
define $b_\xi$ to be the restriction \Res {f(\ASeq \xi)} {(\Prod {\seq
\kappa} \Minus A_\xi)}, where
\[
        A_\xi = \Set [\Big] {\tau \in \Prod {\seq \kappa}} {f(\ASeq
        \xi)(\tau) \in \BigUnion [\zeta < \xi] {\Ran [\big] {f(\ASeq
        \zeta)}}}.
\]
Since $\lambda_\eta \leq \kappa_0$ and for every $\xi \not= \zeta$,
the set \Set [\big] {\tau(0)} {\tau \in \Prod {\seq \kappa} \And
f(\ASeq \xi)(\tau) = f(\ASeq \zeta)(\tau)} is nonstationary in
$\kappa_0$, $A_\xi$ is in \NonstatIdeal {\seq \kappa}. Hence \Seq
{b_\xi} {\xi < \lambda_\eta} is a transversal for \Quotient {\B \eta}
{\A [^*] \eta}, contrary to the choice of \B \eta.

Because the sequences in \B \eta are finite, the standard compactness
argument \Note {for the first order logic} yields a finite subset
$F_\eta$ of \B \eta such that \Quotient {F_\eta} {\A [^*] \eta} does
not have such a choice function
\Note {the neatness is needed because a similar argument cannot be
applied to transversals}.

Now we can define for every $\eta \in \Sf[']$, that $I_\eta$ is such
that $\Set {\ASeq i}
{i \in I_\eta} = F_\eta$ and
\[
        \ASeq \eta' = \Seq [\big] 
        {a_{i,l} \Comp \pi_{\eta,i,l}} {
                i \in I_\eta \And
                a_{i,l} \in \ASeq i \AInter \A [^*] \eta}.
\]

We show first that $\AFamily['] = \Set {\ASeq \eta'} {\eta \in
\Sf[']}$ is almost free. Fix $J$ to be a subset of \Sf['] having
cardinality $< \lambda$. Since the given \AFamily is almost free,
there is a transversal \Seq {b_i} {i \in \BigUnion [\eta \in J]
{I_\eta}} for \Set {\ASeq i} {i \in \BigUnion [\eta \in J] {I_\eta}}.
We define a transversal \Seq {b_\eta'} {\eta \in J} for \Set {\ASeq
\eta'} {\eta \in J}.
Fix an arbitrary \eta from $J$.  Assume, toward a contradiction, that
for every $i \in I_\eta$, $b_i$ is a restriction of the coordinate
$a_{i,l_i}$ in \Quotient {\ASeq i} {\A [^*] \eta}. But then for all $i
\not= j \in I_\eta$, $a_{i,l_i} \not= a_{j,l_j}$ contrary to the
choice of $F_\eta$.
Therefore, there is $i \in I_\eta$ such that for some coordinate
$a_{i,l}$ from the nonempty sequence $\ASeq i \AInter \A [^*] \eta$,
$b_i$ is a restriction of $a_{i,l}$. Thus $X_{\eta, i,l}$ is not in
\NonstatIdeal {\seq \kappa}.
By \Property {PiMap}, $\Prod {\seq \kappa} \Minus \Dom {b_i} \in
\NonstatIdeal {\seq \kappa}$ implies that the set
$\Prod {\seq \kappa} \Minus Y_{\eta,i}$ is in \NonstatIdeal {\seq
\kappa}, where $Y_{\eta,i} = \InvImage {\pi_{\eta,i,l}} {\Dom {b_i}
\Inter X_{\eta, i,l}}$.
Hence we may define $b_\eta' = \Res {b_i \Comp \pi_{\eta,i,l}}
{Y_{\eta,i}}$.

It remains to show that all the other demands listed in \Definition
{NRTSkeleton} can be fulfilled. Exactly as for the \NPT \lambda
{-}-skeletons, \LS['], \LSystem['], and \AFamily['] can be modified so
that \LS['] is in the canonical form of type \seq \lambda and height
\Star n, and moreover, \LSystem['] is disjoint.

By choosing a suitable sub-\lambda-set of \LS['] if necessary \cite
[Claim 3.2(1)] {Sh161}, we may assume that for all $n < \Star n$ and
$\eta, \nu \in \Sf[']$, the sets
\Set [\big] {l < \Lh {\ASeq \eta'}} {\Ran {a_{\eta,l}'} \Subset \LSet
['] {\Res \eta {n+1}}}
and \Set [\big] {l < \Lh {\ASeq \nu'}} {\Ran {a_{\nu,l}'} \Subset
\LSet ['] {\Res \eta {n+1}}}
are equal. Hence there exists \Star l such that $\Lh {\ASeq \eta'} =
\Star l$ for every $\eta \in \Sf[']$, and also, the required partition
\Seq {L^{n+1}} {n < \Star n} of \Star l exists.

When $\Alephf 0 < \theta_n < \min \seq \lambda$ is well-defined one
can add $\theta_n$ into \seq \kappa \Note {see \Fact {AddIndex}
below}, if not yet appeared there. So assume that $\theta_n$ has the
index $k$ in \seq \kappa and define $a_{\eta,l}''(\tau) = \Pair
{a_{\eta,l}'(\tau)} {\alpha_{\tau(k)}}$ for every $l < \Star l$ and
$\tau \in \Prod {\seq \kappa}$, where \Seq {\alpha_\xi} {\xi <
\kappa_k} is an increasing sequence of ordinals cofinal in $\eta(n)$
\Note {change \LSystem['] accordingly}.

For the rest of the properties \Note {covering the case $\theta_n >
\max \seq \kappa$}, a suitable modification procedure is described in
a very detailed way in \cite [Theorem VII.3A.5] {EklofMekler} \Note
{e.g., the case that $\theta_n > \min \seq \lambda$ and $\theta_n
\not\in \seq \lambda$}.
%
%
\end{Proof}
\end{LEMMA}

\begin{COROLLARY}{NRTImpliesNPT}%
\begin{ITEMS}

\ITEM{Implies}%
Suppose \NRTSkeleton is an \NRT \lambda {-}-skeleton of type \seq
\lambda. Then there is a family \SFamily of countable sets such that
the triple
\[
        \SimpleSeq [\Big] {
                \LS, 
                \Seq [\big] {
                        \SubsetsOfCard {\LSet \rho} {< \Alephf 0}
                } {\rho \in \LS \Minus \Sf},
                \SFamily%
        } 
\]
forms an \NPT \lambda {-}-skeleton.

\ITEM{Equal}%
For every uncountable cardinal \lambda, \NPT \lambda - holds if and
only if there is a neat family exemplifying that \NRT \lambda - holds.

\end{ITEMS}

\begin{Proof}%
The reader may wonder why a family witnessing that \NRT \lambda -
holds is not ``straightforwardly'' an example of a family witnessing
that \NPT \lambda - holds. The answer is that, because \NRT -
{-}-transversal is a stronger notion than the standard transversal,
the non-freeness in the \NRT - {-}-sense does not guarantee the
non-freeness in the standard sense.

\ProofOfItem{Implies}%
The existence of an \NPT \lambda {-}-skeleton follows from \Fact
{NotFree} together with \cite [Claim 3.8] {Sh161} or \cite [VII.3A.6]
{EklofMekler}.
To see that the structure of the \lambda-system can be almost
preserved, consider a family \Set {\s \eta} {\eta \in \Sf} such that
each coordinate \Note {\wrt to the slightly modified \lambda-system}
\slevel \eta {n+1} corresponds to a ``well-chosen'' subset of all
finite sequences of elements in \BigUnion [l \in L^{n+1}] {\Ran
{a_{\eta,l}}} \Note {the reason to use finite sequences is the same as
in the proof of \Proposition {Singular}}.

\ProofOfItem{Equal}%
\relax From left to right the claim follows from \cite [Claim 3.8] {Sh161} or
\cite [VII.3A.6] {EklofMekler}.
The other direction follows from \Lemma {NRTSkeleton} and \Item
{Implies}.
\end{Proof}%
\end{COROLLARY}

The following two small facts are needed in the proof of \Lemma
{Transitive}.

\begin{FACT}{TransIffSmall}%
For all \NRT \lambda {\seq \kappa}-skeletons \NRTSkeleton and $I
\Subset \Sf$, $I$ is small in \Sf if and only if there exists an \NRT
[\NonstatIdeal {\seq \kappa}]-{\seq \kappa}-transversal for \Set
{\ASeq \eta} {\eta \in I}.

\begin{Proof}%
A proof of this fact is similar to a proof of the property
\Property [\RefString {Section} {TransGame}] {Trans}
in \Section {TransGame}.
\end{Proof}%
\end{FACT}

\begin{FACT}{AddIndex}
If there is an \NRT \lambda {\seq \kappa}-skeleton of type \seq
\lambda, then for every regular cardinal $\chi < \min \seq \lambda$,
there exists an \NRT \lambda {\seq \chi}-skeleton of type \seq
\lambda, where $\Ran {\seq \chi} = \Ran {\seq \kappa} \Union
\Singleton \chi$.

\begin{Proof}%
Suppose $\seq \chi = \Seq {\chi_k} {k < \Star k}$ is a decreasing
enumeration of $\Ran {\seq \kappa} \Union \Singleton \chi$ and $K$
denotes the index set \Set {k < \Star k} {\chi_k \not= \chi}.
For every $\eta \in \Sf$ and $l < \Lh {\ASeq \eta}$, replace the old
coordinate $a_{\eta,l}$ with a new one $a_{\eta,l}'$, where
$a_{\eta,l}'$ has domain \Prod {\seq \chi} and for every $\tau \in
\Prod {\seq \chi}$, $a_{\eta,l}' (\tau) = a_{\eta,l}(\Res \tau K)$.
Define a new family \AFamily['] to be \Set {\ASeq \eta'} {\eta \in
\Sf}, where $\ASeq \eta'$ is \Seq {a_{\eta,l}'} {l < \Lh {\ASeq
\eta}}.  The only problem is that \AFamily['] should be almost \NRT
[\NonstatIdeal {\seq \chi}]-{\seq \chi}-free.
However, for any $X \Subset \Prod {\seq \kappa}$, $\Prod {\seq \kappa}
\Minus X \in \NonstatIdeal {\seq \kappa}$ implies that \Set {\tau \in
\Prod {\seq \chi}} {\Res \tau K \not\in X} is in \NonstatIdeal {\seq
\chi} as well.
Hence every \NRT [\NonstatIdeal {\seq \kappa}]-{\seq
  \kappa}-transversal \Seq {b_\eta} {\eta \in I} can be
straightforwardly transformed into the form of an \NRT [\NonstatIdeal
{\seq \chi}]-{\seq \chi}-transversal \Seq {b_\eta'} {\eta \in I},
where each \Dom {b_\eta'} equals \Set [\big] {\tau \in \Prod {\seq
    \chi}} {\Res \tau K \in \Dom {b_\eta}}.
\end{Proof}%
\end{FACT}

In the light of the previous fact and for purposes of the forthcoming
``transitivity'' lemma, it makes sense to define ``compatible
skeletons''.

\begin{DEFINITION}{Compatible}
Suppose \NRTSkeleton ['] is an \NRT \lambda {\seq \kappa}-skeleton of
type \seq \lambda and \NRTSkeleton [''] is an \NRT \sigma {\seq
  \chi}-skeleton of type \seq \sigma. We say that \NRTSkeleton ['] is
compatible with \NRTSkeleton [''] if the following conditions are
satisfied:
\begin{textprops}

\TextItem{type}%
$\min \seq \lambda > \max \seq \sigma$;

\TextItem{index}%
for all cardinals $\kappa \in \Ran {\seq \kappa}$, if $\kappa \geq
\min \seq \sigma$ then there is $n$ below the height of \LS[''] such
that for every $\rho \in \LS ['']$ of length $n$, $\sigma_\rho =
\kappa$. 

\end{textprops}
\end{DEFINITION}

\begin{LEMMA}{Transitive}%
Suppose that an \NRT \lambda {\seq \kappa}-skeleton \NRTSkeleton [']
of type \seq \lambda is compatible with an \NRT \sigma {\seq
  \chi}-skeleton \NRTSkeleton [''] of type \seq \sigma.
Then there exists an \NRT \lambda {\seq \kappa \Union \seq
  \chi}-skeleton \NRTSkeleton of type \Concat {\seq \lambda} {\seq
  \sigma}.

\begin{Proof}
%
By \Fact {AddIndex} we may assume that \seq \chi is an end segment of
\seq \kappa \Note {does not have any effect on compatibility and if
\seq \chi is \SimpleSeq {\Alephf 0}, just make sure that \Alephf 0 is
the last element of \seq \kappa}.
By the assumption on \seq \sigma, let $N \Subset n''$ be such that the
initial segment of \seq \kappa of cardinals not in \seq \chi is equal
to \Seq {\sigma_n} {n \in N}.
Thus for every $\eta'' \in \Sf['']$ and $\tau \in \Prod {\seq \chi}$,
the concatenation \Concat {\Res {\eta''} N} \tau is a sequence in
\Prod {\seq \kappa}.

We define \NRTSkeleton in a straightforward manner by
``concatenating'' \NRTSkeleton ['] and several copies of \NRTSkeleton
[''].
So \Sf is defined to be \Set {\Concat {\eta'} {\eta''}} {\eta' \in
\Sf['] \And \eta'' \in \Sf['']}, and \LS is \Set {\Res \eta n} {n <
\Star n}, where $\Star n = n' + n''$, $n'$ is the height of \LS['],
and $n''$ is the height of \LS[''].

To define the other components, fix an arbitrary $\eta = \Concat
{\eta'} {\eta''}$ from \Sf \Note {where this notation means $\eta' \in
\Sf[']$ and $\eta'' \in \Sf['']$}.
Let \Star l denote $l' + l''$, where $l'$ and $l''$ are the lengths of
the sequences in \AFamily['] and \AFamily[''] respectively.

Define $\ASeq \eta = \Seq {a_{\eta,l}} {l < \Star l}$, by setting
for each $l < l'$ and $\tau \in \Prod {\seq \chi}$,
\[
        a_{\eta,l}(\tau)
        =
        \SimpleSeq [\Big] {
                l, \eta'',
                a_{\eta',l}'\Par [\big]{
                        \Concat {\Res {\eta''} N} \tau
                }
        },
\]
and for every $l$ with $l' \leq l < \Star l$ and $\tau \in \Prod {\seq
\chi}$,
\[
        a_{\eta,l}(\tau)
        =
        \SimpleSeq [\big] {
                l, \eta', a_{\eta'',l - l'}''(\tau)
        }.
\]

Define a \lambda-system by setting for every $\rho' \in \LS['] \Minus
\Sf[']$ and $\alpha < \lambda_{\rho'}$,
\[
        \LSet {\ConcatSimple {\rho'} \alpha} =
        l' \times \Sf[''] \times
        \LSet ['] {\ConcatSimple {\rho'} \alpha},
\]
and for every $\eta' \in \Sf[']$, $\rho'' \in \LS[''] \Minus \Sf['']$,
and $\alpha < \sigma_{\rho''}$, set
\[
        \LSet {\ConcatSimple {\Concat {\eta'} {\rho''}} \alpha} =
        \Star l \times
        \Singleton {\eta'}
        \times
        \LSet [''] {\ConcatSimple {\rho''} \alpha}.
\]
Clearly \LS and \LSystem have the desired form and $\AFamily =\Set
{\ASeq \eta} {\eta \in \Sf}$ is based on \LSystem. So to prove that
\NRTSkeleton is an \NRT \lambda \theta-skeleton, it remains to show
that \AFamily is almost \NRT [\NonstatIdeal {\seq \chi}]-{\seq
  \chi}-free.

Fix $I \Subset \Sf$ of cardinality $< \lambda$. We define an
\NRT[\NonstatIdeal {\seq \chi}]-{\seq \chi}-transversal \Seq {b_\eta}
{\eta \in I} for \Set {\ASeq \eta} {\eta \in I}.
The set $I' = \Set {\Res \eta {n'}} {\eta \in I}$ has cardinality $<
\lambda$.
Therefore, there is an \NRT [\NonstatIdeal {\seq \kappa}]-{\seq
  \kappa}-transversal \Seq {b_{\eta'}} {\eta' \in I'} for \Set {\ASeq
  {\eta'}'} {\eta' \in I'}.
Fix now $\eta' \in I'$. Define $J_{\eta'}$ to be the set of all
$\eta'' \in \Sf['']$ such that
\[
        \Set [\big] {\tau \in \Prod {\seq \chi}} {
                \Concat {\Res {\eta''} N} \tau
                \not\in \Dom {b_{\eta'}}
        } \in \NonstatIdeal {\seq \chi}.
\]
Since $\Prod {\seq \kappa} \Minus \Dom {b_{\eta'}}$ is in
\NonstatIdeal {\seq \kappa} and \seq \chi is an end segment of \seq
\kappa, the complement $K_{\eta'} = \Sf[''] \Minus J_{\eta'}$ must be
small in \Sf['']
\Note {remember, $X \Subset \Sf['']$ is small in \Sf[''] if \Set {\Res
{\eta''} N} {\eta'' \in X} is in \NonstatIdeal {\Res {\seq \kappa} k},
where $k$ is the largest index below \Lh {\seq \kappa} with $\kappa_k
> \max \seq \chi$}.
By \Fact {TransIffSmall} there exists an \NRT[\NonstatIdeal {\seq
  \chi}]-{\seq \chi}-transversal \Seq {b_{\eta', \eta''}} {\eta'' \in
  K_{\eta'}} for \Set {\ASeq {\eta''}''} {\eta'' \in K_{\eta'}}.

Consider $\eta \in \Sf$. If $\eta''$ is in $J_{\eta'}$, then define
\[
        b_\eta = \Res {a_{\eta,l}} {
                \Set [\Big] {\tau \in \Prod {\seq \chi}}{
                        \Concat {\Res {\eta''} N} \tau
                        \in \Dom {b_{\eta'}}
                }
        },
\]
where $l < l'$ is the index with $b_{\eta'} \Subset a_{\eta',l}'$.
Otherwise $\eta''$ is in $K_{\eta'}$, and we can define
\[
        b_\eta = \Res {a_{\eta,l}} {\Dom {b_{\eta', \eta''}}}.
\]%
where $l = l' + m$ and $m < l''$ is the index with $b_{\eta', \eta''}
\Subset a_{\eta'',l}''$.%
%
%
%
%
\end{Proof}

\end{LEMMA}

\begin{COROLLARY}{Transitive}
For all regular cardinals $\lambda > \kappa > \chi$, if \NRT \lambda
\kappa and \NRT \kappa \chi hold, then \NRT \lambda \chi holds.
\end{COROLLARY}
%

\end{SECTION}

\begin{SECTION}{Conclusion}{Conclusion}

%
Now we may put the pieces together. By \Definition
{CanonicalLambdaSet}, if \lambda is a regular cardinal below the
possible first regular limit cardinal and \LS is a \lambda-set in the
canonical form of height \Star n and type \Seq {\lambda_\rho} {\rho
\in \LS \Minus \Sf},
then there exist sequences $\seq \lambda = \Seq {\lambda_n} {n < \Star
n}$ and $\seq \theta = \Seq {\theta_n} {n < \Star n}$ such that for
every $\rho \in \LS \Minus \Sf$ of length $n$,
both $\lambda_\rho = \lambda_n$ and $\LSE \rho \Subset \Cof
{\theta_n}$ hold. For the rest of this section we assume that all the
types of skeletons are given in this simplified form.

\begin{REMARK}
For simplicity we have chosen \Prod {\seq \kappa} to be the
domain of the functions appearing in the elements of \NRT-{\seq
\kappa}-families. Hence we must restrict ourselves below the possible
first regular limit cardinal.  It is possible to replace \Prod {\seq
\kappa} with a ``full'' $\kappa_0$-set and \NonstatIdeal {\seq \kappa}
with ``small sets'' \wrt the fixed $\kappa_0$-set, but that makes the
notation unnecessarily complicated.
\end{REMARK}

\begin{DEFINITION}{NRTSet}%
%
%
Let \NRTSet denote the smallest set of cardinals such that
\begin{properties}

\item \Alephf 0 is in \NRTSet;

\item if there exists an \NRT \lambda {\seq \kappa}-skeleton of type
\seq \lambda such that there is no regular limit cardinal below
\SuccCard \lambda and $\Ran {\seq \kappa} \Subset \NRTSet$,
then $\Ran {\seq \lambda} \Subset \NRTSet$.
\end{properties}
\end{DEFINITION}

\begin{REMARK}
One could expect that, analogously to \Fact {NPTImpliesNRT}
and because of the transitivity property \Lemma {Transitive}, it would
suffice to consider \NRT \lambda {\seq \kappa}-skeletons of height 1
only \Note {i.e., that in the definition above we could assume \seq
  \lambda has length 1}. However, an analogous proof does not work for
\NRT \lambda {\seq \kappa}-skeletons because, even in a fixed ``level
$n$'', a transversal may choose from several possible coordinates
$a_{\eta,l}$, $l \in L^{n+1}$, contrary to the \NPT \lambda {-}-case,
where on each level $n$ there is only one coordinate, namely \slevel
\eta {n+1}.
\end{REMARK}

The benefit of \Fact {NPTImpliesNRT} is that we can separate levels of
\NPT \lambda {-}-skeletons to independent building blocks and combine
those blocks in various ways. Of course in this countable case, we
apply \Corollary {NRTImpliesNPT} \Note {i.e., the coordinates
  $a_{\eta,l}$, $l \in L^{n+1}$, of level $n$ in the definition of an
  \NRT-{-}-family can be amalgamated to a single coordinate}.

\begin{LEMMA}{AddLevel}%
Suppose $\sigma < \lambda$ are cardinals in \NRTSet. There exists an
\NPT \lambda {-}-skeleton of type \seq \lambda such that if \sigma is
uncountable, $\sigma \in \Ran {\seq \lambda}$.

\begin{Proof}
If \lambda is the first uncountable cardinal in \NRTSet, i.e.,
$\lambda = \Alephf 1$, then the claim holds by the definition. Suppose
\NRTSkeleton[^0] is an \NRT \lambda {\seq \kappa}-skeleton of type
\seq \chi, $\Ran {\seq \chi} \Union \Ran {\seq \kappa} \Subset
\NRTSet$, and the claim holds for all cardinals below \lambda. Note
that the problematic case is $\max \seq \chi > \sigma > \min \seq
\chi$.

By the induction hypothesis there is an \NPT {\max \seq \kappa}
{-}-skeleton \NPTSkeleton[^1] of type \seq \theta such that $\Ran
{\seq \kappa} \Subset \Ran {\seq \theta}$ and if $\max \seq \kappa >
\sigma > \Alephf 0$ then $\sigma \in \Ran {\seq \theta}$ holds too.
Since \NRTSkeleton[^0] and \NPTSkeleton[^1] are compatible, it follows
from \Lemma {Transitive} and \Corollary {NRTImpliesNPT}, that there is
an \NPT \lambda {-}-skeleton \NPTSkeleton[^2] of type \Concat {\seq
  \chi} {\seq \theta}.

Suppose $\sigma > \min \seq \chi$ \Note {hence $\sigma \not\in \Ran
{\seq \theta}$} and $n$ is the largest index with $\chi_n > \sigma >
\chi_{n+1}$.
By the induction hypothesis there is an \NPT \sigma {-}-skeleton
\NPTSkeleton[^3] of such a type that it contains the sequence
$\seq \sigma = \Concat {\Concat {\SimpleSeq \sigma} {\Seq {\chi_m} {n
< m < \Star n}}} {\seq \theta}$.
By applying \Fact {NPTImpliesNRT} and \Fact {AddIndex} to
\NPTSkeleton[^2], there is an \NRT \lambda {\seq \sigma}-skeleton
\NRTSkeleton[^4] of type \Seq {\chi_m} {m \leq n}.
Now \NPTSkeleton[^3] and \NRTSkeleton[^4] are compatible, and hence
the claim follows from \Lemma {Transitive} together with \Corollary
{NRTImpliesNPT}.
\end{Proof}
\end{LEMMA}


\begin{CONCLUSION}{NicelyIncompact}%
\begin{ITEMS}

\ITEM{NPT}%
For all cardinals \lambda in \NRTSet, \NPT \lambda {-} hold. In
particular, \NRTSet is a nicely incompact set of regular cardinals in
the sense of \cite [Definition 2.1] {Sh523}.

\ITEM{Sh523}%
If $K$ is a nicely incompact set of regular cardinals below the
possible first inaccessible cardinal, then $K \Subset \NRTSet$.

\ITEM{MgSh204}%
All the cardinals from the set \MgShSet defined in \cite [Theorem 1 of
\SS1] {MgSh204} belong to \NRTSet \Note {or look at \Theorem
  {SimplifiedForm}}.
Particularly, all the regular cardinals below $\Alephf {\omega \Times
\omega +1}$ are in \NRTSet, and \NRTSet is cofinal below the first
cardinal fixed point.

\end{ITEMS}

\begin{Proof}%
\ProofOfItem{NPT}%
By \Lemma {AddLevel}.

\ProofOfItem{Sh523}%
Suppose \seq \chi is a type of some \NPT \lambda {-}-skeleton. By
induction on increasing order of the cardinals in \seq \chi, apply
\Fact {NPTImpliesNRT} to show that $\Ran {\seq \chi} \Subset \NRTSet$.

\ProofOfItem{MgSh204}%
By \ItemOfFact {ExistenceOfNRT} {Successor} and \ItemOfFact
{ExistenceOfNRT} {MgSh204}.
\end{Proof}%
\end{CONCLUSION}

Our final conclusion concerns game-free groups. As we have seen in
\Proposition {Singular}, game-freeness of an \NPT \lambda {-}-skeleton
\NPTSkeleton is closely connected to the possible value of \LSE
\emptyset.
Hence we have to look at a little bit more restricted set of nicely
incompact cardinals.

\begin{DEFINITION}{NRTGameFreeSet}%
Let \NRTGameFreeSet denote the set of all cardinals \lambda in \NRTSet
such that \lambda appears in a type of some \NRT \lambda {\seq
\kappa}-skeleton \NRTSkeleton satisfying that $\seq \kappa \Subset
\NRTGameFreeSet$,
and moreover, if \lambda is a successor of a singular cardinal \mu,
then $\LSE \emptyset \Inter \Cof {\Cf \mu} = \emptyset$.
\end{DEFINITION}

\begin{FACT}{Difference}%
All the cardinals from the set \MgShSet \Note {\ItemOfConclusion
  {NicelyIncompact} {MgSh204}} belong to \NRTGameFreeSet.
In fact, all the examples in \Fact {ExistenceOfNRT}, except the first
two of them, yield skeletons fulfilling the cofinality demand in the
definition of \NRTGameFreeSet.

\end{FACT}

\begin{THEOREM}{game-free}%
Suppose \mu is a cardinal such that both \Cf \mu and $\lambda =
\SuccCard \mu$ are in \NRTGameFreeSet.
Then there exists a nonfree $(\mu \Times \theta)$-game-free group of
cardinality \lambda, where $\theta < \mu$ is a regular cardinal such
that if \mu is a successor of a regular cardinal then $\mu = \SuccCard
\theta$, and otherwise $\theta \not= \Cf \mu$.

\begin{Proof}%
The proof proceeds by induction on increasing order of the cardinals
in \NRTGameFreeSet.
By \cite [\AuxRef {PV1} {Lemma} {GameFreeGroup}] {PV1}, \cite [\AuxRef
{PV1} {Lemma} {SingularCompactness}] {PV1}, and \cite [\AuxRef {PV1}
{Lemma} {RegForm}] {PV1} it suffices to show existence of \NPT \lambda
{-}-skeleton such that if \mu is regular then \SFamily is
\mu-game-free, and if \mu is singular, then \SFamily is
\epsilon-game-free for every $\epsilon < \mu$.

The case successor of a regular cardinal follows from \cite [\AuxRef
{PV1} {Fact} {SColFree}\ and \AuxRef {PV1} {Lemma} {RegForm}] {PV1},
since as an induction hypothesis, we may assume that there exists a
\sigma-game-free \NPT \mu {-}-skeleton \NPTSkeleton, where \sigma is a
cardinal such that $\mu = \SuccCard \sigma$ and if \sigma is a regular
cardinal, then $\sigma = \SuccCard \theta$ and $\LSE \emptyset \Subset
\Cof \theta$.

The case successor of a singular cardinal follows from \Conclusion
{NicelyIncompact} and \Proposition {Singular}, since $\NRTGameFreeSet
\Subset \NRTSet$ and a modification of \Lemma {AddLevel} for
\NRTGameFreeSet holds \Note {thus the demands of \Proposition
  {Singular} can be fulfilled}.
\end{Proof}%
\end{THEOREM}
%

\end{SECTION}

\begin{SECTION}{CubGame}{On cub-game and game-free groups}

%
In this section \mu is a singular cardinal, \kappa denotes the
cofinality of \mu, and \lambda is the successor cardinal of \mu.
We study existence of nonfree \epsilon-game-free families of
cardinality \lambda, when $\epsilon > (\mu \Times \mu) + \mu$.
Our tools for that are some known and modified results about
``cub-game'' for successors of singular cardinals. An analog study for
successors of regular cardinals, which is an easier case, is carried
out in \cite [\AuxRef {PV1} {Section} {Basic}] {PV1}.

\begin{DEFINITION}{CubGame}%
Suppose \lambda is an uncountable regular cardinal $A \Subset
\lambda$, and $\epsilon < \lambda$ is an ordinal. For notational purposes
let $x$ denote a set of regular cardinals below \lambda \Note {$x$
tells in which ``cofinalities the limits are checked''}.
We denote by \CubGame [x] \epsilon A \lambda the following two players
cub-game. The players, \PlayerI \Note {also called ``outward'' player}
and \PlayerII \Note {also called ``inward'' player}, choose in turns a
sequence \Seq {\alpha_i} {i < \epsilon} of ordinals such that
\begin{properties}

\item \PlayerI chooses ordinals $\alpha_0 < \lambda$ and $\alpha_{i+2}
< \lambda$, for $i < \epsilon$, with $\alpha_{i+2} > \alpha_{i+1}$;

\item if $i$ is a limit ordinal, \PlayerI must choose the ordinal
$\alpha_i = \sup_{j < i} \alpha_j$;

\item when $\alpha_i$ is defined, \PlayerII must pick some ordinal
$\alpha_{i+1} < \lambda$ with $\alpha_{i+1} > \alpha_i$.

\end{properties}
\PlayerII [big] wins a play if for all limit ordinals i for which $\Cf
i \in x$, $\alpha_i$ belongs to $A$.

We say that \PlayerII wins \CubGame [x] \epsilon A \lambda if there exists
a winning strategy for \PlayerII in \CubGame [x] \epsilon A \lambda.
For a regular cardinal \theta below \lambda, \CubGame [\theta] \epsilon A
\lambda is a shorthand for \CubGame [\Singleton \theta] \epsilon A
\lambda.
\end{DEFINITION}

\begin{DEFINITION}{good}%
\AddInTheorem{%
\cite [Def. 1.9] {Sh523}, %
\cite [Def. 2 of \SS1] {MgSh204}, %
\cite [\SS 1] {Sh351}%
}
Suppose \Pair {\seq \mu} {\scale f} is a scale for a singular cardinal
\mu of cofinality \kappa \Note {\Definition {scale}}, where $\seq \mu
= \Seq {\mu_\xi} {\xi < \kappa}$ and $\scale f = \Seq {f_\alpha}
{\alpha < \lambda}$.
We denote by \good {\scale f} the set of all ``good points \alpha \wrt
\scale f'', i.e., \alpha's below \lambda such that for some unbounded
set $A \Subset \alpha$ and $\zeta < \kappa$, $f_\alpha(\xi) <
f_\beta(\xi)$ holds whenever $\alpha < \beta \in A$ and $\zeta \leq
\xi < \kappa$.
Let \bad {\scale f} denote the set of all limit ordinals $\alpha <
\lambda$ such that $\alpha \not\in \good {\scale f}$ and $\Cf \alpha >
\kappa = \Cf \mu$.
\end{DEFINITION}

\begin{PROPOSITION}{ExtremelyGameFreeGroups}
Assume that the following conditions are fulfilled:
\begin{properties}

\item \mu is a singular cardinal, $\lambda = \SuccCard \mu$, and $\Cf
\mu = \kappa$;

\item there is a scale \Pair {\seq \mu} {\scale f} for \mu such that
\bad {\scale f} is stationary in \lambda;

\item there exists an \NPT \lambda {-}-skeleton \NPTSkeleton of type
\seq \lambda such that $\LSE \emptyset \Subset \bad {\scale f}$;

\item if \kappa is uncountable then there is $n$ below the height of
\LS such that for every $\rho \in \LS$ of length $n$, $\lambda_\rho =
\kappa$ \Note {where $\lambda_\rho$ is from the type \seq \lambda}.

\end{properties}
Then there exists a nonfree group of cardinality \lambda, which is
\delta-game-free for every $\delta < \lambda$
\Note {and even in a ``longer'' game}.

\begin{Proof}%
By \Proposition {Singular} we may assume that \NPTSkeleton is chosen
so that \SFamily is \epsilon-game-free for every $\epsilon < \mu$.
By \cite [\AuxRef {PV1} {Lemma} {SingularCompactness}] {PV1} we know
that \SFamily is \mu-game-free.
Fix a \delta below \lambda. As in \cite [\AuxRef {PV1} {Lemma}
{RegForm}] {PV1}, we can fix a filtration \Seq {\S \alpha} {\alpha <
  \lambda} of \SFamily, and ensure using a suitable bookkeeping, that
after any ``block'' of \mu moves by the players of the transversal
game \TransGame \delta \SFamily, there is some $\alpha < \lambda$ such
that the elements chosen by the players are from \S \alpha, and
moreover, all the elements of \S \alpha has been chosen.

Suppose \theta is a regular cardinal below \mu such that $\LSE
\emptyset \Subset \Cof \theta$. Note that $\LSE \emptyset \Subset \bad
{\scale f}$ implies $\theta > \kappa$.  By \Lemma {IIwinsGood} and
\Lemma {LongCubGame}, \PlayerII has a winning strategy in the cub-game
\CubGame [\theta] \delta {\good{\scale f}} \lambda.

During the transversal game \TransGame \delta \SFamily, \PlayerII can
additionally use his winning strategy in the cub-game \CubGame
[\theta] \delta {\good{\scale f}} \lambda to ensure that after
arbitrary many rounds of the blocks of \mu moves in the transversal
game, the elements chosen by the players are exactly the elements in
\S \alpha for some \alpha which is not in $\LSE \emptyset = \Set
{\beta < \lambda} {\Quotient{\SFamily}{\S \beta} \Text {is not
\lambda-free}} \Subset \bad {\scale f} \Inter \Cof \theta$.

By \cite [\AuxRef {PV1} {Lemma} {Quotient}] {PV1}, the family
$\Quotient {\SFamily} {\S \alpha} = \Set {s \Minus \Union{\S \alpha}}
{s \in \SFamily \Minus \S \alpha}$ is \mu-game-free. Therefore
\PlayerII can continue the transversal game \TransGame \delta \SFamily
one more round of the block of \mu moves. During these new \mu moves
\PlayerII uses his bookkeeping and winning strategy in the cub-game
again, and so on, up to all the required \delta moves.
Now the claim follows from \cite [\AuxRef {PV1} {Lemma}
{GameFreeGroup}] {PV1}.
\end{Proof}%
\end{PROPOSITION}

Before changing the subject to the winning strategies in the cub-game,
we ask: for which singular cardinals \mu can the demands of the last
proposition be fulfilled? We need few lemmas before the conclusion.

\begin{LEMMA}{supercomp->bad_stat}%
\AddInTheoremCite[Claim 27]{Sh108}%
Suppose \chi is a supercompact cardinal, \mu is a singular cardinal,
$\kappa = \Cf \mu$, $\kappa < \chi < \mu$, $\lambda = \SuccCard \mu$,
and \Pair {\seq \mu} {\scale f} is a scale for \mu.
There exists a singular cardinal $\sigma < \mu$ such that
$\Cf \sigma = \kappa$ and
$\bad {\scale f} \Inter \Cof {\SuccCard \sigma}$ is stationary in
\lambda.

\begin{Proof}%
Let us first show that \bad {\scale f} is stationary in \lambda.
Suppose, contrary to this claim, that \C is a cub of \lambda with $\C
\Inter \bad {\scale f} = \emptyset$ and $j$ is an embedding from
\GroundModel onto an inner model $M$ satisfying that $j(\xi) = \xi$
for every $\xi < \chi$, $j(\chi) \geq \lambda$, and $\SubsetsOfCard M
{\leq \lambda} \Subset M$.
Let \delta be $\sup \Image j \lambda$. The set \Image j \lambda is in
$M$ and \delta is in $j(\C)$, because \C is a cub of \lambda, $\delta
< \sup j(\C) = j(\lambda)$, and $\Image j \lambda \Subset j(\C)$.
So \delta is good \wrt $j(\scale f)$ for $j(\mu)$ in $M$ and \Image j
\lambda is a cofinal subset of \delta. By \cite [Lemma 6] {MgSh204}
there exists a cofinal subset $A$ of \Image j \lambda and \zeta
witnessing the goodness of \delta.
Define $A'$ to be the set \InvImage j A. Then for every $\alpha <
\beta \in A'$, the inequality $f_\alpha(\zeta) < f_\beta(\zeta)$
holds, since $j$ is an elementary embedding, $j(\zeta) = \zeta$,
$j(\alpha), j(\beta) \in A$, and ${(j(f))}_{j(\alpha)}(\zeta) <
{(j(f))}_{j(\beta)}(\zeta)$.
This is a contradiction, since $A'$ has cardinality \lambda, and the
set \Set [\big] {f_\alpha(\zeta)} {\alpha < \lambda} has cardinality
$\mu_\zeta < \mu < \lambda$ \Note {where $\mu_\zeta$ is from \seq
\mu}.

Now the desired conclusion follows from the fact that \lambda is a
successor of the singular cardinal \mu of cofinality \kappa and \delta
has cofinality \lambda in $M$: If $\alpha \in \C$ whenever $\sup \C
\Inter \alpha = \alpha$ and \alpha has cofinality \rho, where \rho is
a successor of a singular cardinal having cofinality \kappa, then
\Image j \C has this property too. So the assumption $\C \Inter \bad
{\scale f} = \emptyset$ leads to a contradiction as above.
\end{Proof}%
\end{LEMMA}

\begin{LEMMA}{BadIsNonreflecting}
Suppose \kappa is a regular cardinal, $\mu = \Alephf \kappa > \kappa$,
$\lambda = \SuccCard \mu$, \Pair {\seq \mu} {\scale f} is a scale for
\mu, and moreover, $2^\kappa = \SuccCard \kappa$.
Then $\bad {\scale f} \Subset \Cof {\SuccCard \kappa}$ and $\bad
{\scale f} \Inter \alpha$ is nonstationary in \alpha, for every
$\alpha < \lambda$.

\begin{Proof}%
Suppose \alpha below \lambda has cofinality $\sigma > \SuccCard
\kappa$. We show that \alpha is good \Note {of course, all the
ordinals of cofinality $< \kappa$ are good}.
By \cite [II.1.2A(3)] {CardArith}, $2^\kappa = \SuccCard \kappa$
implies that there exists an exact upper bound $g$ of \Seq {f_\beta}
{\beta < \alpha} \Note {i.e., $g$ is an \EventLess-upper bound for
\Seq {f_\beta} {\beta < \alpha} such that for every $h \in \Prod {\seq
\mu}$, $h \EventLess g$ implies that $h \EventLess f_\beta$ for some
$\beta < \alpha$}.

Now argue as in \cite [Case 1 of the proof of Lemma 5 in \SS1]
{MgSh204}: By \cite [Lemma 7] {MgSh204} the set \Set [\big] {\xi <
\kappa} {\Cf [\big] {g(\xi)} > \sigma} is in \BoundedIdeal \kappa.
Assume, toward a contradiction, that \alpha is not good. By \cite
[Lemma 6] {MgSh204} \Seq [\big] {\Cf [\big] {g(\xi)}} {\xi < \kappa}
is not eventually constant, even though, its range is a subset of
$\Alephf \kappa \Inter (\sigma +1)$ modulo \BoundedIdeal \kappa.
Since there are $< \kappa$ cardinals in the range of this sequence,
there must exists different cardinals $\sigma_1$ and $\sigma_2$ such
that both $\Set [\big] {\xi < \kappa} {\Cf [\big] {g(\xi)} = \sigma_1}
\not\in \BoundedIdeal \kappa$ and $\Set [\big] {\xi < \kappa} {\Cf
[\big] {g(\xi)} = \sigma_2} \not\in \BoundedIdeal \kappa$ hold.
However, by \cite [Lemma 8] {MgSh204}, $\sigma_1 = \sigma = \sigma_2$
must hold, a contradiction.

The claim on nonreflection follows from the definition of goodness: if
\alpha is in \good {\scale f}, then there is a closed unbounded subset
of \alpha consisting of ordinals in \good {\scale f} only.
\end{Proof}%
\end{LEMMA}

Using \cite [Fact 4.2] {Sh351} \Note {similarly to \cite [Conclusion
29] {Sh108}} we get the following conclusion.

\begin{LEMMA}{Collapsing}
Suppose \chi is a supercompact cardinal, \GCH holds, $\kappa < \Alephf
\kappa$ is a regular cardinal below \chi, \mu is \SuccCard [\kappa]
\chi, and \lambda is \SuccCard \mu.
Then there is a forcing extension, where $\ZFC + \GCH$ holds, no
bounded subset of \kappa is added, \mu is the singular cardinal
$\Alephf \kappa$, $\lambda = \SuccCard \mu$, and all the assumptions
of \Proposition {ExtremelyGameFreeGroups} hold too.

\begin{Proof}%
If \NPT \kappa {-} does not hold, shoot a nonreflecting stationary
subset $F$ of $\kappa \Inter \Cof {\Alephf 0}$ by a forcing notion
described, e.g., in \cite [Proof of Lemma 3.1] {Sh521}. This is for
building the desired \NPT \lambda {-}-skeleton at the end of this
proof. Then no bounded subset of \kappa is added, \chi remains
supercompact, and \GCH still holds.

By \Lemma {supercomp->bad_stat}, there is a cardinal \theta which is a
successor of a singular cardinal $\sigma < \lambda$ so that $\Cf
\sigma = \kappa$ and $\bad {\scale f} \Inter \Cof \theta$ is
stationary in \lambda. Let $E_1$ denote this stationary set. Note that
by the form of \sigma and \chi, $\sigma < \SuccCard \sigma = \theta <
\chi$ holds.

Using Levy collapse \Col \kappa {< \sigma}, collapse all the cardinals
below \sigma to \kappa.
Because of \GCH and $\Cf \sigma = \kappa$,
\Col \kappa {<\sigma} has cardinality \sigma in \GroundModel. Recall
that \Col \kappa {<\sigma} is \kappa-complete. So in \VCollapse \kappa
{< \sigma}, no bounded subset of \kappa is added, \kappa is a regular
cardinal, $F$ is a nonreflecting stationary subset of \kappa, $E_1$ is
still a stationary subset of \lambda, \sigma has cardinality \kappa,
\theta is the cardinal \SuccCard \kappa, and \GCH holds.
Moreover, the \kappa-completeness and $\Card [\big] {\Col \kappa
{<\sigma}} = \sigma < \theta < \lambda$ implies that
\Pair {\seq \mu} {\scale f} is still a scale for \mu,
and as mentioned in \cite [4.2(2)] {Sh351}, if we let $E_2$ denote the
set \bad {\scale f} in \VCollapse \kappa {< \sigma}, $E_1 = E_2$
modulo a cub of \lambda in \VCollapse \kappa {< \sigma}.

Now use the Levy collapse \Col \theta {<\chi}, to collapse all the
cardinals between $\theta = \SuccCard \kappa$ and \chi. Since in
\VCollapse \kappa {<\sigma}, \chi is still a strongly inaccessible
cardinal, \Col \theta {<\chi} has \chi-c.c.
Consequently, in the final extension, $E_2$ is still a stationary
subset of \lambda, $\chi = \SuccCard \theta = \SuccCard [2] \kappa$,
$\mu = \Alephf \kappa$ \Note {by the assumption $\kappa < \Alephf
\kappa$}, and $\lambda = \SuccCard \mu$.
Since \Col \theta {<\chi} is also \theta-complete, it follows that in
the final extension, no bounded subset of \kappa is added, \Pair {\seq
\mu} {\scale f} remains as a scale for \mu, \GCH holds, $\theta =
\SuccCard \kappa$, $F$ is still a nonreflecting stationary subset of
\kappa, and by \cite [4.2(2)] {Sh351} again, $\bad {\scale f} = E_2$
modulo a cub of \lambda. Consequently \bad {\scale f} is stationary in
$\lambda = \Alephf {\kappa +1}$.

It remains to show that the required type of \NPT \lambda {-}-skeleton
exists.
By the properties of $F$, \ItemOfFact {ExistenceOfNRT} {Nonreflecting}
and \ItemOfFact {ExistenceOfNRT} {Successor}, we get an \NPT \theta
{-}-skeleton in whose type \kappa appears, if \kappa is uncountable.
By \Lemma {BadIsNonreflecting}, \ItemOfFact {ExistenceOfNRT}
{Nonreflecting}, and \Lemma {Transitive} we get an \NPT \lambda
{-}-skeleton \NPTSkeleton in whose type \theta appears, and also
\kappa is in the type if \kappa is uncountable.
Furthermore, the set \LSE \emptyset equals $\lambda \Inter \Cof
\theta$. Therefore, by shrinking \LSE \emptyset to \bad {\scale f'},
one gets the required \NPT \lambda {-}-skeleton.
\end{Proof}%
\end{LEMMA}


In ``extreme cases'' \good {\scale f} does not contain a closed and
unbounded subset of \lambda \Note {for more cases, see \cite [Fact
1.7] {Sh351} and \cite [Claim 4.3] {FoMa1997}}.
However, \good {\scale f} is ``almost'' a cub of \lambda, because
there is a winning strategy for the ``inward'' player in a very long
cub-game.

\begin{LEMMA}{IIwinsGood}%
Suppose $\lambda = \SuccCard \mu$ and $\mu > \Cf \mu = \kappa$. We
denote by \CubGame [\not=\kappa] \epsilon A \lambda the game \CubGame
[x] \epsilon A \lambda, where $x$ is the set of all other regular
cardinal below \lambda except \kappa.
For every $\epsilon < \mu$, \PlayerII has a winning strategy in the game
\CubGame [\not=\kappa] {\epsilon} {\good{\scale f}} \lambda.

\begin{Proof}%
First of all let \Star \zeta be the least index below \kappa with
$\epsilon < \mu_{\Star \zeta}$.
Suppose that the players has already chosen the ordinals \Seq
{\alpha_j} {j \leq i} and \PlayerII should choose $\alpha_{i+1}$.
Assume that \PlayerII has picked during the earlier rounds also
functions \Seq {h_j} {j < i \Text{is odd}} satisfying for each odd $j$
that
\begin{properties}

\item $h_j \in \Prod [\xi < \kappa] \mu_\xi$;

\item for every odd $k < j$ and $\xi \in \kappa \Minus \Star \zeta$,
both $f_{\alpha_k}(\xi) < h_j(\xi)$ and $h_k(\xi) < h_j(\xi)$ hold;

\item $h_j \EventLess f_{\alpha_j}$.

\end{properties}
Firstly, \PlayerII defines $h_{i+1}$ by setting $h_{i+1}(\xi) = 0$ for
all $\xi < \Star \zeta$, and otherwise,
\[
        h_{i+1}(\xi) =
        \Par[\Big]{
                \sup_{j < i \Text{ odd}}
                \Par[\big]{\Max{f_{\alpha_j}(\xi), h_j(\xi)}}
        } +1.
\]
Then $h_{i+1}$ is in \Prod [\xi < \kappa] {\mu_\xi}, since $\mu_\xi$
is a regular cardinal greater than \epsilon when $\xi \geq \Star \zeta$.
Secondly, \PlayerII picks $\alpha_{i+1}$ to be the least ordinal \beta
above $\alpha_i$ satisfying that $h_{i+1} \EventLess f_\beta$. Such an
ordinal exists because \scale f is cofinal in \Prod [\xi < \kappa]
{\mu_\xi}.

So it remains to show that for a limit $i$ of cofinality not equal to
\kappa, $\alpha_i = \Sup [j<i] {\alpha_j}$ belongs to \good {\scale
f}.
If $\Cf i < \kappa$, $\alpha_i$ is good. So assume that $\Cf i >
\kappa$.
Let $I$ be a cofinal subset of $i$ having order type \Cf i and
consisting of odd ordinals only \Note {the moves of \PlayerII}.
For every $j \in I$, define $\xi_j < \kappa$ to be the smallest index
with $h_j \StrictlyLess[\xi_j] f_{\alpha_j}$. Since $\Cf i > \kappa$,
there is a cofinal subset $J$ of $I$ and \zeta such that for all $j
\in J$, $\xi_j = \zeta$. But then $A = \Set {\alpha_j} {j \in J}$ is a
cofinal subset of $\alpha_i$ and $A$ together with \zeta witness that
$\alpha_i$ is good, since for all $k < j$ from $J$ and for all $\xi
\geq \zeta$,
\[%
        f_{\alpha_k}(\xi)
        < h_j(\xi)
        < f_{\alpha_j}(\xi).%
\]%
\end{Proof}%
\end{LEMMA}

\begin{REMARK}
If $2^\mu = \lambda$ then the conclusion of the previous lemma follows
from the definition of \IGoodIdeal \lambda too, as explained in
detail, e.g., in \cite [Lemma 2.1] {HuHyRa2000}.
\end{REMARK}

\begin{LEMMA}{LongCubGame}%
Assume $\kappa \leq \mu$ are cardinals, \lambda is \SuccCard \mu, $A
\Subset \lambda$, and $x$ a set of regular cardinals below \lambda.
Suppose that \PlayerII wins \CubGame [x] \epsilon A \lambda for every
$\epsilon < \mu$, and moreover, if \mu is a regular cardinal,
\PlayerII wins \CubGame [x] {\mu+1} A \lambda.
Then \PlayerII wins \CubGame [x] \sigma A \lambda for every $\sigma <
\lambda$.

\begin{Proof}%
This is a known fact, presented e.g. in \cite [\AuxRef {PV1} {Section}
{Basic}] {PV1}.
\end{Proof}%
\end{LEMMA}

It is again a known fact that the inward player has a winning strategy
even in ``much longer'' games of the form \CubGame [x] T A \lambda,
where the length of the game is measured by the tree $T$ of closed
subsets of $A$ of order type $\alpha+1$, $\alpha < \lambda$ \Note
{ordered by the end extension}.
Details are presented in \cite [\AuxRef {PV1} {Section} {Basic}]
{PV1}.

%

\end{SECTION}



%
\begin{tabbing}
Saharon Shelah:\= \\
        \>Institute of Mathematics\\
        \>The Hebrew University\\
        \>Jerusalem, Israel\\
\\
        \>Department of Mathematics\\
        \>Rutgers University\\
        \>New Brunswick\\
        \>NJ, USA\\
        \>\texttt{shelah@math.rutgers.edu}
\end{tabbing}

\begin{tabbing}
Pauli V\"{a}is\"{a}nen:\= \\
        \>Department of Mathematics\\
        \>University of Helsinki\\
        \>Finland\\
        \>\texttt{pauli.vaisanen@helsinki.fi}
\end{tabbing}
%

\end{document}